\documentclass[11pt]{article}
\usepackage{amsmath,amssymb}
\usepackage[latin1]{inputenc}
\usepackage[dvips]{graphicx}

\topskip=\baselineskip  \headsep=0in
\newtheorem{theorem}{Theorem}[section]
\newtheorem{lemma}[theorem]{Lemma}
\newtheorem{proposition}[theorem]{Proposition}

\newtheorem{corollary}[theorem]{Corollary}
\newtheorem{remark}[theorem]{Remark}

\newcommand{\qed}{\enspace\vrule height6pt width4pt depth2pt}

\newenvironment{proof}{\par\noindent{\bf Proof.}}{$\qed$\par\bigskip}

\newcommand{\vcd}{{\rm vcd}}

\newcommand{\Cen}{\mbox{\rm C}}
\newcommand{\core}{\mbox{\rm Core}}
\newcommand{\Aut}{\mbox{\rm Aut}}

\newcommand{\inv}{^{-1}}

\newcommand{\Z}{{\mathbb Z}}
\newcommand{\Q}{{\mathbb Q}}
\newcommand{\R}{{\mathbb R}}
\newcommand{\C}{{\mathbb C}}
\newcommand{\HQ}{{\mathbb H}}

\newcommand{\SL}{{\rm SL}}

\newcommand{\PSL}{{\rm PSL}}

\newcommand{\matriz}[1]{\begin{array} #1 \end{array}}

\newcommand{\GEN}[1]{\langle #1 \rangle}

\newcommand{\gorro}[1]{\widehat{#1}}

\newcommand{\quat}[2]{\left( \frac{#1}{#2} \right)}

\newcommand{\U}{{\cal U}}
\newcommand{\V}{{\cal V}}
\newcommand{\T}{{\cal T}}
\newcommand{\W}{{\cal W}}
\newcommand{\CC}{{\cal C}}

\title{Groups of units of integral group rings
commensurable with direct products of free-by-free
groups \footnote{Research partially supported by
Onderzoeksraad of Vrije Universiteit Brussel, Fonds
voor Wetenschappelijk Onderzoek (Belgium), D.G.I. of
Spain, Fundaci{\'{o}}n S{\'{e}}neca of Murcia, CNPq
and FAPDF. \newline AMS classification index:
Primary 16U60, Secondary 11R27, 16A26.}}
\author{Eric Jespers \and Antonio Pita \and {\'{A}}ngel del R\'{i}o
\and Manuel Ruiz \and Pavel Zalesski}

\date{}
\begin{document}

\def\thedefinition{\null}
\def\thetheorem{\arabic{theorem}}

\maketitle

\begin{abstract}
We classify the finite groups $G$ such that the
group of units of the integral group ring $\Z G$ has
a subgroup of finite index which is a direct product
of free-by-free groups.
\end{abstract}

The investigations on the unit group $\Z G^{*}$ of
the integral group ring $\Z G$ of a finite group $G$
have a long history and go back to work of Higman
\cite{HigmaThesis}. One of the fundamental problems
that attracts a lot of attention is Research Problem
17 posed by Sehgal in \cite{seh-book2}: find
presentations of $\Z G^{*}$ for some finite groups
$G$. For many finite groups $G$ a finite set of
generic generators of a subgroup of finite index in
$\Z G^{*}$ has been obtained, but there is no
general result known on determining the relations
among these generators. This work was initiated by
Bass and Milnor \cite{bass-milnor} and then
continued by Kleinert \cite{K}, Ritter and Sehgal
\cite{RS}, and Jespers and Leal \cite{JL}. For a
survey on the above mentioned results, we refer to
Sehgal's book \cite{seh-book2} and to
\cite{jespers-sur}.

An alternative approach to that of finding presentations is the
one suggested by Kleinert in \cite{KleinertSurvey}. Recall that a
generic unit group of $A$ is a subgroup of finite index in the
group of reduced norm $1$ elements of an order in $A$. Then
according to Kleinert ``a unit theorem for  a finite dimensional
simple rational algebra $A$ consists of the definition, in purely
group theoretical terms, of a class of groups $\CC(A)$ such that
almost all generic unit groups of $A$ are members of $\CC(A)$''.
This approach has an obvious generalization to finite dimensional
semi-simple rational algebras, such as the rational group algebra
$\Q G$ of a finite group $G$ and its orders, for example $\Z G$.
This kind of unit theorem has been obtained for integral group
rings $\Z G$ of some restricted classes of finite groups $G$.  We
give a brief history on the descriptions obtained so far. Higman
in \cite{HigmaThesis} showed that if $G$ is a finite abelian group
then $\Z G^*=L\times (\pm G)$, where $L$ is a free abelian group
of rank depending on the cardinality of $G$ and the order of the
elements of $G$. This result heavily depends on Dirichlet's Unit
Theorem. He also showed that if $G$ is non-abelian then $\Z G^*$
is finite if and only if $G$ is a Hamiltonian 2-group and in this
case $\Z G^*=\pm G$. The finite groups $G$ such that $\Z G^*$ is
virtually free and non-abelian (there are only four) were
classified in \cite{J}. This last result was motivated by a
previous theorem of Hartley and Pickel \cite{HP} which states that
$\Z G^*$ is either abelian, finite or has a non-abelian free
subgroup. Finally, the finite groups $G$ such that $\Z G^*$ is
virtually a direct product of free groups (there are infinitely
many) were classified in a series of papers by Jespers, Leal and
del R{\'\i}o \cite{JLR,JR,LR}. Thus the finite groups $G$ for
which a unit theorem, in the sense of Kleinert,  is known for $\Z
G^*$ are those for which the class of groups considered are either
finite groups, abelian groups, free groups or direct products of
free groups. As far as we know, all the finite groups $G$ for
which the structure of $\Z G^*$ is known up to commensurability
are covered by these results.

The aim of this paper is to obtain a group
theoretical description of $\Z G^*$ for a larger
family of finite groups $G$ than the family of
groups mentioned in the previous paragraph. We do
this  by connecting the study of $\Z G^*$ with the
better known structure of the Bianchi groups. The
inspiration came from some examples in
\cite{corrales} and \cite{PRR}. In the second
reference some presentations of $\Z G^*$ are
obtained for two groups of order $16$ for which $\Z
G^*$ is commensurable with the Picard group
$\PSL_2(\Z[i])$. These two groups belong to a class
of finite groups, called groups of Kleinian type,
for which geometrical methods are applicable to
obtain presentations of groups of finite index
(implementation of the method however is usually
difficult). Our main theorem (Theorem~\ref{Main})
shows that the class $\CC$ containing the generic
groups of $\Z G^*$ for $G$ of Kleinian type is
formed by the direct products of free-by-free
groups, and in fact this property characterizes the
groups of Kleinian type. Furthermore, we classify
the finite groups of Kleinian type as the groups
which are epimorphic images of some specific groups.
This classification is the most involving part of
the paper. In order to state this result we first
fix some terminology.

\vspace{12pt} Recall that a group $H$ is said to be
{\em free-by-free} if $H$ contains a normal subgroup
$N$ so that both $N$ and $H/N$ are free groups. Note
that the trivial and infinite cyclic group are free
groups, and thus free groups and finitely generated
abelian groups are direct products of free-by-free
groups.

For a ring $R$ we denote by $R^*$ the group of
invertible elements of $R$ and by $Z(R)$ its centre.
In case $R$ is an order in a simple finite
dimensional rational algebra $A$ we denote by $R^1$
the group consisting of the elements of reduced norm
$1$ in $R$. (By an order we always mean a
$\Z$-order; see \cite{seh-book2} for a definition.)

Two subgroups $H_{1}$ and $H_{2}$ of a group $H$
are said to be {\em commensurable} when their
intersection has finite index in both $H_{1}$ and
$H_{2}$. Often the group $H$ is clear from the
context and hence will not be specifically
mentioned. For instance, the statement ``$\Z G^*$
is commensurable with a direct product of
free-by-free groups'' means that $\Z G^*$ is
commensurable with some subgroup of $\Q G^*$ which
decomposes in a direct product of free-by-free
groups. Similarly, if $R$ is an order in a simple
finite dimensional rational algebra $A$, then the
statement ``$R^{1}$ is commensurable with a
free-by-free group'' means that $R^{1}$ is
commensurable with a subgroup of $A^{*}$ with the
mentioned property.

A finite group $G$ is said to be of {\em Kleinian
type} if every non-commutative simple quotient $A$
of the rational group algebra $\Q G$ has an
embedding $\psi:A\rightarrow M_2(\C)$ such that
$\psi(R^1)$ is a discrete subgroup of $\SL_2(\C)$
for some (every) order $R$ in $A$.

It turns out that if $G$ is a finite group of
Kleinian type then it is metabelian. Hence to state
our classification of the groups of Kleinian type it
is convenient to introduce some notation for
presentations of such groups. The cyclic group of
order $n$ is usually denoted by $C_n$. To emphasize
that $a\in C_n$ is a generator of the group, we
write $C_n$ either as $\GEN{a}$ or $\GEN{a}_n$.
Recall that a group $G$ is metabelian if $G$ has an
abelian normal subgroup $N$ such that $A=G/N$ is
abelian. We simply denote this information as
$G=N:A$. To give a concrete presentation of $G$ we
will write $N$ and $A$ as direct products of cyclic
groups and give the necessary extra information on
the relations between these generators. By
$\overline{x}$ we denote the coset $xN$. For
example, the dihedral group of order $2n$ and the
quaternion group of order $4n$ can be described as
    $$\matriz{{llll}
    D_{2n} & = & \GEN{a}_n : \GEN{\overline{b}}_2, & b^2 = 1, a^b = a\inv. \\
    Q_{4n} & = & \GEN{a}_{2n} : \GEN{\overline{b}}_2, & a^b = a\inv, \; b^2 = a^n.}$$
If $N$ has a complement in $G$ then $A$ can be
identify with this complement and we write
$G=N\rtimes A$. For example, the dihedral group also
can be given by $D_{2n} = \GEN{a}_n \rtimes
\GEN{b}_2$ with $a^b = a\inv$.

We are now in a position to formulate the main
result.

    \topmargin=-0.4in 
    \textheight=9.5in 

\newpage



\begin{theorem}\label{Main}
For a finite group $G$ the following statements are
equivalent.
\begin{itemize}
\item[(A)] $\Z G^*$ is
commensurable with a direct product of free-by-free
groups.
\item[(B)]
For every simple quotient $A$ of $\Q G$ and some
(every) order $R$ in $A$, $R^1$ is  commensurable
with a free-by-free group.
\item[(C)]
For every simple quotient $A$ of $\Q G$ and some
(every) order $R$ in $A$, $R^1$ has virtual
cohomological dimension at most 2.
\item[(D)]
$G$ is of Kleinian type.
\item[(E)]
Every simple quotient of $\Q G$ is either a field, a
totally definite quaternion algebra or $M_2(K)$,
where $K$ is either $\Q$, $\Q(i)$, $\Q(\sqrt{-2})$
or $\Q(\sqrt{-3})$.
\item[(F)]
$G$ is either abelian or an epimorphic image of
$A\times H$, where $A$ is abelian and one of the
following conditions holds:
\begin{enumerate}
\item $A$ has exponent $6$ and $H$ is one of the following groups:

$\bullet$
 $\W = \left(\GEN{t}_2 \times \GEN{x^2}_2
\times \GEN{y^2}_2\right) :
\left(\GEN{\overline{x}}_2 \times
\GEN{\overline{y}}_2\right)$,
 with $t = \left(y,x\right)$ and $Z(\W)=\GEN{x^2,y^2,t}$.

$\bullet$
 $\W_{1n} = \left(\prod\limits_{i=1}^n
\GEN{t_i}_2 \times \prod\limits_{i=1}^n \GEN{y_i}_2
\right) \rtimes \GEN{x}_4$,
 with $t_{i} = (y_i,x)$ and $Z(\W_{1n})=\GEN{t_1,\dots,t_n,x^2}$.

$\bullet$
 $\W_{2n} = \left(\prod\limits_{i=1}^n
\GEN{y_i}_4 \right) \rtimes \GEN{x}_4$, with $t_{i}
= (y_i,x) = y_i^2$ and
$Z(\W_{2n})=\GEN{t_1,\dots,t_n,x^2}$.

\item $A$ has exponent $4$ and $H$ is one of the following groups:

$\bullet$ $\V = \left(\GEN{t}_2 \times \GEN{x^2}_4
\times \GEN{y^2}_4\right) :
 \left(\GEN{\overline{x}}_2 \times \GEN{\overline{y}}_2\right)$,
 with $t = (y,x)$ and $Z(\W)=\GEN{x^2,y^2,t}$.

$\bullet$
 $\V_{1n} = \left(\prod\limits_{i=1}^n
\GEN{t_i}_2 \times \prod\limits_{i=1}^n \GEN{y_i}_4
\right) \rtimes \GEN{x}_8$,
 with $t_{i} = (y_i,x)$ and $Z(\V_{1n})=\GEN{t_1,\dots,t_n,y_1^2,\dots,y_n,x^2}$.

$\bullet$ $\V_{2n} = \left(\prod\limits_{i=1}^n
\GEN{y_i}_8 \right) \rtimes \GEN{x}_8$,
 with $t_{i} = (y_i,x) = y_i^4$ and $Z(\V_{2n})=\GEN{t_i,x^2}$.

$\bullet$ $\U_1 = \left(\prod\limits_{1\le i< j \le
3} \GEN{t_{ij}}_2 \times \prod\limits_{ k=1}^3
\GEN{y_k^2}_2\right) :\left(\prod\limits_{k=1}^3
\GEN{\overline{y_k}}_2\right)$,
 with $t_{ij} = (y_j,y_i)$ and $Z(\U_1)=\GEN{t_{12},t_{13},t_{23},y_1^2,y_2^2,y_3^2}$

$\bullet$ $\U_2 = \left(\GEN{t_{23}}_2 \times
\GEN{y_1^2}_2 \times \GEN{y_2^2}_4 \times
\GEN{y_3^2}_4 \right) :
 \left(\prod\limits_{k=1}^3 \GEN{\overline{y_k}}_2 \right)$,
 with $t_{ij} = (y_j,y_i)$, $y_2^4=t_{12}$, $y_3^4=t_{13}$ and
 $Z(\U_2)=\GEN{t_{12},t_{13},t_{23},y_1^2,y_2^2,y_3^2}$.

\item $A$ has exponent $2$ and $H$ is one of the following groups:

$\bullet$
 $\T= \left(\GEN{t}_4 \times
\GEN{y}_8\right) : \GEN{\overline{x}}_2$, with $t =
(y,x)$ and $x^2 = t^2 = (x,t)$.

$\bullet$
 $\T_{1n} = \left(\prod\limits_{i=1}^n
\GEN{t_i}_4 \times \prod\limits_{i=1}^n \GEN{y_i}_4
\right) \rtimes \GEN{x}_8$,
 with $t_{i} = (y_i,x)$, $(t_i,x)=t_i^2$ and $Z(\T_{1n})=\GEN{t_1^2,\ldots,t_n^2,x^2}$.

$\bullet$
 $\T_{2n} = \left(\prod\limits_{i=1}^n \GEN{y_i}_8 \right) \rtimes \GEN{x}_4$,
 with $t_{i} = (y_i,x) = y_i^{-2}$ and $Z(\T_{2n})=\GEN{t_1^2,\ldots,t_n^2,x^2}$.

$\bullet$
 $\T_{3n} = \left( \prod\limits_{i=2}^n
 \GEN{t_i}_4 \times
 \GEN{y_1^2t_1}_2 \times \GEN{y_1}_8\times
 \prod\limits_{i=2}^n \GEN{y_i}_4 \right) :
 \GEN{\overline{x}}_2$, with $t_{i} = (y_i,x)$,
 $(t_i,x)=t_i^2$, $x^2 = t_1^2$, and
 $Z(\T_{3n})=\GEN{t_1^2,\ldots,t_n^2,x^2}$.

 \item
 $H = M\rtimes P = (M\times
 Q):\GEN{\overline{u}}_2$, where $M$ is an
 elementary abelian $3$-group,
 $P=Q:\GEN{\overline{u}}_2$, $m^u=m\inv$ for every
 $m\in M$, and one of the following conditions
 holds:

 $\bullet$ $A$ has exponent $4$ and $P=C_8$.

 $\bullet$ $A$ has exponent $6$, $P=\W_{1n}$ and
 $Q=\GEN{y_1,\dots,y_n,t_1,\dots,t_n,x^2}$.

 $\bullet$ $A$ has exponent $2$, $P=\W_{21}$ and $Q=\GEN{y_1^2,x}$.
\end{enumerate}

\end{itemize}

\end{theorem}

\renewcommand{\baselinestretch}{1}\large\mbox{}\normalsize
    \topmargin=-0.15in 
    \textheight=9in 


According to \cite{SerreGal}, a group $G$ is called
{\it good} if the homomorphism of cohomology groups
\newline $H^n(\widehat G,M)\longrightarrow H^n( G,M)$
induced by the natural homomorphism
$G\longrightarrow \widehat G$ of $G$ to its
profinite completion $\widehat G$ is an isomorphism
for every finite $G$-module $M$.

Theorem~\ref{Main} yields that for a finite group
$G$  of Kleinian type  the non-commutative simple
components of $\Q G$ that are not totally definite
quaternion algebras are of the form $M_{2}(\Q
(\sqrt{-d}))$ with $d=0,1,2$ or $3$. On the other
hand  the groups of units of an order in a number
field and in a totally definite quaternion algebra
are commensurable with a free abelian group.
Therefore, since the group of units of two orders
in $\Q G$ are commensurable, Theorem~\ref{Main}
implies that $\Z G^*$ is commensurable with a
direct product of a free abelian group and groups
of the form $\SL_2(\Z[\sqrt{-d}])$ with $d=0,1,2$
or $3$. Recently it was shown that
$\PSL_2(\Z[\sqrt{-d}])$ is good for every
non-negative integer \cite{GJZ}. Since the class of
good groups is closed under commensurability and
$\PSL_2(\Z[\sqrt{-d}])$ has a subgroup of finite
index isomorphic to a subgroup of finite index of
$\SL_2(\Z[\sqrt{-d}])$, it follows that
$\SL_2(\Z[\sqrt{-d}])$ is good. Moreover the class
of good groups is closed under finite direct
products. Hence the following property follows at
once.

\begin{corollary}
If $G$ is a finite group of Kleinian type then the
group of units of its integral group ring $\Z G$ is
good.
\end{corollary}

In particular, this corollary says that the  virtual
cohomological dimension of the profinite completion
of $ZG^*$ coincides with the virtual cohomological
dimension of $ZG^*$ and so the profinite completion
of $ZG^*$ is virtually torsion free.

\medskip

The outline of the paper is as follows. In Section 1
we introduce the basic notation used throughout the
paper. In Section 2 we show that conditions (A) and
(B) are equivalent. In Section 3 we prove (B)
implies (C) (which is obvious), (C) implies (D) (by
first classifying the simple algebras of Kleinian
type and the finite dimensional simple algebras $A$
for which $R^1$ has virtual cohomological dimension
at most 2 for an order $R$ in $A$) and (E) implies
(B) (by using known facts about Euclidean Bianchi
groups). Section 4 is dedicated to prove (F) implies
(E). At this point one has shown that all the groups
satisfying condition (F) are of Kleinian type. The
most involved part of the proof is to show that (D)
implies (F), that is showing that condition (F)
exhausts the class of groups of Kleinian type. This
is proved for nilpotent groups in Section 5 and for
non-nilpotent groups in Section 6.

In a preliminary version of the proof of (D) implies
(F) we used previous results from \cite{JR,PRR}. We
thank Jairo Gon\c{c}alves for attracting our
attention to an old result of Amitsur which
classifies the finite groups that have all its
irreducible complex characters of degree $1$ or $2$.
This result has been very helpful in reducing
earlier given arguments and in making the proof of
(D) implies (F) independent of \cite{JR,PRR}.

\section{Preliminaries}

\def\thetheorem{\thesection.\arabic{theorem}}

We introduce the basic notation and the main tools
used in the paper. The Euler function is denoted by
$\varphi$. For a positive integer $n$, let $\xi_n$
denote a complex primitive root of unity.

Let $G$ be a group.  For $x,y\in G$, we put
$x^{y}=y^{-1}xy$ and $(x,y)=xyx^{-1}y^{-1}$. We
recall the following well known formulas:
$(ab,c)=(b,c)^{a\inv} (a,c)$ and
$(a,bc)=(a,b)(a,c)^{b\inv}$. The centre and derived
subgroup of $G$ are denoted by $Z(G)$ and $G'$
respectively. The notation $H\le G$ means that $H$
is a subgroup of $G$ and if $H$ is a normal subgroup
of $G$ then we write $H\unlhd G$. The normalizer of
$H\le G$ in $G$ is denoted by $N_G(H)$. If $N\unlhd
G$ then we will use the usual bar notation for the
natural images of the elements and subsets of $G$ in
$G/N$, that is $\overline{x}$ denotes the coset $xN$
of $x\in G$ and if $X\subseteq G$ then
$\overline{X}$ denotes $\{\overline{x}\mid x\in
X\}$. A semidirect product associated to an action
of a group $H$ on a group $N$ is denoted by
$N\rtimes H$.

We say that a group  {\em virtually satisfies} a
group theoretical condition if it has a subgroup of
finite index satisfying the given condition. For
example, $G$ is virtually abelian if and only if
$G$ has an abelian subgroup of finite index. Notice
that if a class of groups satisfying a property
$\mathcal{P}$ is closed under subgroups of finite
index then a group $G$ is commensurable with a
group satisfying $\mathcal{P}$ if and only if it
virtually satisfies $\mathcal{P}$.  Moreover, in
this case, $G$ is commensurable with a group which
is a direct product of groups satisfying
$\mathcal{P}$ if and only if it is virtually a
direct product of groups satisfying $\mathcal{P}$.
This, of course, applies to the class of
free-by-free groups.

As well as the groups described in statement (F) of
Theorem 1, the following metabelian groups will be
relevant. $$\matriz{{llll}
    D_{2^{n+2}}^{+} &=& \GEN{a}_{2^{n+1}} \rtimes \GEN{b}_2, &  a^b = a^{2^n + 1}.  \\
    D_{2^{n+2}}^{-} &=& \GEN{a}_{2^{n+1}} \rtimes \GEN{b}_2,  &  a^b = a^{2^n - 1}. \\
    \mathcal{D}     &=& \left(\GEN{c}_4 \times \GEN{a}_2\right) \rtimes \GEN{b}_2, &
                    Z(\mathcal{D}) = \GEN{c}, \; \left(b,a\right)=c^2. \\
    \mathcal{D}^{+} &=& \left(\GEN{c}_4 \times \GEN{a}_4\right) \rtimes \GEN{b}_2, &
                    Z(\mathcal{D}^{+}) = \GEN{c}, \;(b,a)=ca^2.
    }$$
Recall that if $G$ is a non-abelian group of order
$2^{n+2}$ having a cyclic subgroup of index $2$ then
$G$ is isomorphic to either the dihedral group
$D_{2^{n+2}}$, the quaternion group $Q_{2^{n+2}}$ or
one of the two semi-dihedral groups $D_{2^{n+2}}^+$
or $D_{2^{n+2}}^-$.

If $K$ is a field and $a$ and $b$ are two non zero
elements of $K$ then $\quat{a,b}{K}$ denotes the
quaternion $K$-algebra defined by $a$ and $b$, that
is, the $K$-algebra given by the following
presentation:
    $$\quat{a,b}{K}=K[i,j|i^2=a,j^2=b,ji=-ij].$$
In case $a=b=-1$, then the previous algebra is also
denoted $\HQ(K)$. It is well known that
$\quat{a,b}{K}$ is split, that is, it is isomorphic
to $M_2(K)$, if and only if the equation
$aX^2+bY^2=Z^2$ has a solution different from
$X=Y=Z=0$.

Let $A$ be a finite dimensional semi-simple rational algebra and
$R$ an order in $A$. Then $R^*$ is commensurable with the group of
units of every order in $A$ (see for example
\cite[Lemma~4.6]{seh-book2}). Assume that, furthermore, $A$ is
simple and let $K$ be the centre of $A$. Then $R^*$ is
commensurable with $Z(R)^*\times R^1$, where $R^1$ denotes the
group of elements of reduced norm 1 of $R$. Moreover $\R
\otimes_{\Q} K \cong \R^r \oplus \C^s$, where $r$ is the number of
embeddings of $K$ in $\R$ and $s$ is the number of pairs of non
real embeddings of $K$ in $\C$. These embeddings in $\R$ and pairs
of embeddings in $\C$ correspond to infinite places of $K$ (i.e.
equivalence classes of archimedean valuations of $K$). If $d$ is
the degree of $A$ then $\R\otimes_{\Q} A \cong M_d(\R)^{r_1}
\oplus M_{d/2}(\HQ(\R))^{r_2} \oplus M_d(\C)^s$, where $r_2$ is
the number of infinite places at which $A$ is ramified and
$r=r_1+r_2$. Every embedding $\sigma$ of $K$ in $\C$ induces an
embedding $\bar{\sigma}:A\rightarrow M_d(\C)$ that maps $R^1$ into
$\SL_d(\C)$.

A  totally definite quaternion algebra is a
quaternion algebra $A$ over a totally real number
field $K$ which is ramified at every infinite place,
that is, $\sigma(K)\otimes_K A \cong \HQ(\R)$ for
every embedding $\sigma:K\rightarrow \R$; or
equivalently $A=\quat{a,b}{K}$ with
$\sigma(a),\sigma(b)<0$ for every field homomorphism
$\sigma:K\rightarrow \R$.

The simple algebra $A$ is said to be of {\em
Kleinian type} if there is an embedding
$\psi:A\rightarrow M_2(\C)$ such that $\psi(R^1)$ is
a discrete subgroup of $\SL_2(\C)$ for some (any)
order $R$ in $A$, or equivalently if $A$ is either a
number field or $A$ is a quaternion algebra and
$\bar{\sigma}(R^1)$ is a discrete subgroup of
$\SL_2(\C)$ for some embedding of $K$ in $\C$. More
generally, an {\em algebra of Kleinian type}
\cite{PRR} is by definition a direct sum of simple
algebras of Kleinian type. So, a {\em finite group}
$G$ is of {\em Kleinian type} if and only if the
rational group algebra $\Q G$ is of Kleinian type.

\section{Equivalence of (A) and  (B)}

The equivalence between (A) and (B) is a direct
consequence of the following more general theorem.

\begin{theorem}\label{Components}
Let $A=\prod_{i=1}^n A_i$ be a finite dimensional
rational algebra such that $A_i$ is simple for every
$i$. Let $R$ be an order in $A$ and for every $i$
let $R_i$ be an order in $A_i$. Then $R^*$ is
virtually a direct product of free-by-free groups if
and only if $R_i^1$ is virtually free-by-free  for
every $i$.
\end{theorem}

A group $G$ is said to be {\em virtually
indecomposable} if every subgroup of finite index of
$G$ is indecomposable as a direct product of two
infinite groups. (Note that the terminology should
not be confused with ``having an indecomposable
subgroup of finite index''.)

To prove Theorem~\ref{Components} we need the
following lemma.

\begin{lemma}\label{ffz1}
If $C$ is a free-by-free group which is virtually
indecomposable and not virtually abelian then
$Z(C)=1$.
\end{lemma}

\begin{proof}
Suppose $C$ is a free-by-free group. Then we may
write $C=N\rtimes F$, with $N$ and $F$ free groups.
We first prove that $Z(C)\subseteq N$. Suppose the
contrary. Then it follows that $Z(F)\ne 1$ and thus
$F$ is cyclic. Therefore $\GEN{Z(C),N}$ has finite
index in $C$. As $C$ is not virtually abelian, $N$
is non-abelian, hence $\GEN{Z(C),N} = Z(C)\times N$,
contradicting the virtual indecomposability of $C$.
So, indeed, $Z(C)\subseteq N$. If $Z(C)\neq 1$ then
$N$ is cyclic and $Z(C)\times F$ has finite index in
$C$, again a contradiction.
\end{proof}

\noindent \textbf{ Proof of
Theorem~\ref{Components}}. Since $R$ and
$\prod_{i=1}^n R_i$ are two orders in $A$ and
$R_i^*$ is commensurable with $Z(R_i)^*\times
R_i^1$ for each $i$, one has that $R^*$ and
$\prod_{i=1}^n Z(R_i)^*\times R_i^1$ are
commensurable. The sufficiency of the conditions is
now clear.

Conversely, assume that the direct product
$T=\prod_{x\in X} T_x$ is a subgroup of finite index
of $R^{*}$, where every $T_x$ is a non trivial
free-by-free group. Since the virtual cohomological
dimension of $R^*$ is finite, $X$ is finite and we
can assume without loss of generality that every
$T_x$ is virtually indecomposable and either $T_x$
is cyclic or is not virtually abelian. For every
$x\in X$ let $\pi_x:T\rightarrow T_x$ denote the
projection and let $Y=\{y\in X \mid T_y \mbox{ is
not abelian}\}$. For every $i$ let $S_i=R_i^1\cap T$
and $Z_i=Z(R_i)^*\cap T$. Then $S_i$ is a torsion
free subgroup of finite index in $R_i^1$, $Z_i$ is a
torsion-free subgroup of finite index in $Z(R_i)^*$,
$S_i\cap Z(R_i)^*=1$ and $S=\prod_i Z_i \times S_i$
is a subgroup of finite index in $T$, because
$R_i^1\cap Z(R_i)$ is finite, $\GEN{Z(R_i)^*,R_i^1}$
has finite index in $R_i^*$ and $T$ is a
torsion-free subgroup of finite index in $R^*$.

We claim that if $\pi_z(S_j)$ is not abelian (and
hence infinite) and $H= \left(\prod_{i\neq j}
S_i\right)\times \left(\prod_i Z_i\right)$ then
$\pi_z(H)=1$. Indeed, $C=\pi_z(S)$ is a subgroup of
finite index in $T_z$ and therefore $C$ satisfies
the hypothesis of Lemma~\ref{ffz1}. Thus
$\pi_z(S_j)\cap \pi_z(H)\subseteq Z(C)=1$, because
$(S_j,H)=1$. Then $C=\pi_z(S_j)\times \pi_z(H)$ and
from the virtual indecomposability of $C$ one
deduces that $\pi_z(H)=1$. This finishes the proof
of the claim.

We have to show that each $R_i^1$ is virtually free-by-free or
equivalently that so is $S_i$. By \cite[Theorem 1]{KR}, $S_i$ is
virtually indecomposable. So either $S_i$ is virtually cyclic, and
we are done, or $S_i$ is non-abelian. Assume that $S_i$ is
non-abelian. Hence there is $y\in Y$ such that $\pi_y(S_i)$ is
non-abelian. Assume that $x\in X_i=\{x\in X \mid \pi_x(S_i)\ne
1\}$. Then $\pi_x(S_i)$ has finite index in $T_x$, for otherwise
$T_x$ is not cyclic and so there is at least one $j$ such that
$\pi_x(S_j)$ is non-abelian that gives, by the claim,
$\pi_x(S_i)=1$, a contradiction. Therefore $S_i$ is a subgroup of
finite index in $\prod_{x\in X_i} T_x$. As $S_i$ is virtually
indecomposable, $|X_i|=1$ and therefore $S_i$ is virtually
free-by-free as wanted. $\qed$

\section{(B) implies (C), (C) implies (D), and (E) implies (B)}

It is well known that the virtual cohomological
dimension of a free-by-free group is at most 2 and
so (B) implies (C) is obvious.

(E) implies (B) is a direct consequence of the
following lemma in which we collect known or
recently established facts on the structure of the
group of reduced norm one elements of an order in
some simple rational algebras.

\begin{lemma}\label{GoodBad}
Let $A$ be a simple finite dimensional rational
algebra and $R$ an order in $A$.
\begin{enumerate}
\item $R^1$ is finite if and only if $A$ is a field or a totally
definite quaternion algebra.
\item $R^1$ is virtually free non-abelian if and only if $A=M_2(\Q)$.
\item If $A=M_2(\Q(\sqrt{-d}))$ with $d=1,2,3,7$ or $11$ then
$R^1$ is commensurable with a free-by-free group.
\end{enumerate}
\end{lemma}

\begin{proof}
See e.g. \cite[Lemma 21.3]{seh-book2} for 1 and
\cite{KleinertSurvey} for 2.

3. Let $O_d$ be the ring of integers of $\Q(\sqrt{-d})$. Then
$R^1$ is commensurable with $\SL_2(O_d)$. So it is enough to show
that $\SL_2(O_d)$ is virtually free-by-free. This is well known
for $d=3$, because $\PSL_3(O_d)$ has a subgroup of index $12$
isomorphic to the figure eight knot group, a free-by-infinite
cyclic group (see for example \cite[page 137]{MR}). That
$\SL_2(O_d)$ is virtually free-by-free for $d=1,2,7$ or $11$ has
been  proved in Lemmas 4.2 and 4.3 of \cite{WZ}.
\end{proof}

In order to prove (C) implies (D) we classify in
Proposition~\ref{ClasiKT} the simple algebras of
Kleinian type (correcting Proposition~3.1 in
\cite{PRR} where one possibility was missed by an
error in the proof) and in Proposition~\ref{vcd2} we
classify the simple algebras $A$ for which $R^1$ has
virtual cohomological dimension at most 2 for an
order $R$ in $A$. Then (C) implies (D) follows at
once from these two propositions.

\begin{proposition}\label{ClasiKT}
A finite dimensional rational simple algebra $A$ is
of Kleinian type if and only if it is either a
number field or a quaternion algebra which is not
ramified at at most one infinite place.


In particular, if $A$ is non-commutative and of
Kleinian type then the centre $K$ of $A$ has at most
one pair of complex non-real embeddings and hence
the order of every primitive root of unity in $K$ is
a divisor of $4$ or $6$.
\end{proposition}

\begin{proof}
Let $R$ be an order in $A$ and $K$ the centre of
$A$. Assume first that $A$ is either a field or a
quaternion algebra which is not ramified at at most
one infinite place. If $A$ is a field or a totally
definite quaternion algebra then $R^1$ is finite by
Lemma~\ref{GoodBad} and so $A$ is of Kleinian type.
If $K$ is totally real then $A$ is of Kleinian type
by a theorem of Borel and Harish-Chandra \cite{BH}
(see \cite{Mac}). Otherwise $K$ has exactly one pair
of complex embeddings and $A$ is ramified at all the
real embeddings of $K$. Thus $A$ is of Kleinian type
by \cite[Theorem 10.1.2]{EGM}.

Conversely, assume that $A$ is of Kleinian type.
Then $A$ is either a number field or a quaternion
algebra. In the remainder of the proof we assume
that $A$ is a quaternion algebra.

Let $\sigma_1,\ldots,\sigma_n$ be the set of
representatives up to conjugation of the embeddings
of $K$ in $\C$. Each $\sigma_i$ gives rise to an
embedding $\overline{\sigma_i}:A\rightarrow A_i$
where $A_i=M_2(\C)$ if $\sigma_i$ is not real,
$A_i=M_2(\R)$ if $\sigma_i$ is real and not ramified
and $A_i=\HQ(\R)$ otherwise. We consider $A_i$
embedded in $M_2(\C)$ in the obvious way. Then
$\sigma_i(R^1)\subseteq SL_2(\C)$. Let $R$ be an
order in $A$. Then $\overline{\sigma_i}(R^1)$ is a
discrete subgroup of $\SL_2(\C)$ for some $i$,
because by assumption $A$ is of Kleinian type. We
may assume that $i=1$. Assume that $\sigma_l$ is
either a non real embedding or a non ramified real
embedding and let $f:A\rightarrow \prod_{j\ne k}
A_j$ be the map given by
$f(x)=(\overline{\sigma_j}(x))_{j\ne k}$. Then, by
the Strong Approximation Theorem (see \cite[Theorem
7.12]{PlR} or \cite[Theorem 4.3]{V}), $f(R^1)$ is
dense in $\prod_{j\ne k} A_j^1$ and therefore $k=1$.
This shows that $A$ ramifies at at least $n-1$
places. Hence the result follows.
\end{proof}

\begin{proposition}\label{vcd2}
Let $A$ be a simple finite dimensional rational
algebra and $R$ an order in $A$. Let $\vcd(R^1)$
denote the virtual cohomological dimension of $R^1$.
The following conditions hold.
\begin{enumerate}
\item $\vcd(R^1)=0$ if and only if $A$ is a field or a totally definite quaternion algebra.
\item $\vcd(R^1)=1$ if and only if $A=M_2(\Q)$.
\item $\vcd(R^1)=2$ if and only if $A=M_2(K)$ with $K$ an imaginary quadratic extension of the rationals or
$A$ is a quaternion algebra over a totally real
number field which is not ramified at exactly one
infinite place.
\end{enumerate}
\end{proposition}

\begin{proof}
Let $K=Z(A)$, $r$ the number of embeddings of $K$ in
$\R$, $s$ the number of of non-real embeddings of
$K$ in $\C$, $r_1$ the number of real embeddings of
$K$ at which $A$ is ramified and $r_2=r-r_1$ the
number of real embeddings of $K$ at which $A$ is not
ramified. Set $A=M_n(D)$ where $D$ is a division
ring of degree $d$. Notice that if $d$ is odd then
$r_1=0$.

The sufficiency of the respective conditions easily
can be checked using the following formula for the
virtual cohomological dimension of $R^1$ that can be
deduced from the formulae on pages 220 and 222 and
in Theorem 4 of \cite{KleinertSurvey}:
    \begin{equation}\label{vcdf1}
    \matriz{{rcl}
    \vcd(R^1) & = & r_2\frac{(nd+2)(nd-1)}{2} + r_1\frac{(nd-2)(nd+1)}{2} + s(n^2d^2-1) -n +1\\
    & = & r_2nd + r\frac{(nd-2)(nd+1)}{2} + s(n^2d^2-1)  - n +1.}
    \end{equation}

Conversely, assume that $A$ is not a field and
$\vcd(R^1)\le 2$. By (\ref{vcdf1}) one has
    \begin{equation}\label{vcdf2}
    r_2nd + r\frac{(nd-2)(nd+1)}{2} + s(n^2d^2-1)\le n+1.
    \end{equation}
Since $A$ is not a field, $nd\ge 2$ and therefore
the three summands on the left hand side of
(\ref{vcdf2}) are non-negative, which implies that
each summand at most $n+1$. Hence, since $nd+1\ge
n+1$, we get that  $s(nd-1)(nd+1) = s(n^2d^2-1)\le
n+1$ and thus it follows that either $s=0$ or
$d=s=1$ and $n=2$. In the latter case $r_1=0$ and
since $s(n^2d^2-1)=n+1$, one has that $r_2nd=0$ so
that $r_2=0$. Thus $A=M_2(K)$ where $K$ is an
imaginary quadratic extension of $\Q$.

Assume now that $s=0$, that is, $K$ is totally real.
Now we use $r_2nd\le n+1$ to deduce that either (a)
$r_2=0$, (b) $r_2=d=1$ or (c) $n=r_2=1$ and $d=2$.
We deal with each case separately.

(a) If $r_2=0$ then $r=r_1\ne 0$, that is $A$ is
ramified at every infinite place of $K$. This
implies that $d$ is even. Furthermore
    $$\frac{(nd-2)(nd+1)}{2}\le r\frac{(nd-2)(nd+1)}{2} \le n+1$$
and therefore $(nd-2)(nd+1)\le 2n+2$. Thus
$nd(nd-1)\le 2n+4$ and so $n(d(nd-1)-2)\le 4$. If
$n\ge 2$ then $n(d(nd-1)-2)\ge 2(2\cdot 3 -2)=8$.
Thus $n=1$, that is $A=D$ is a division ring.
Further $d(d-1)\le 6$ and thus $d=2$, because $d$ is
even. We conclude that $A$ is a totally definite
quaternion algebra.

(b) Assume that $r_2=d=1$. Then $r_1=0$, that is,
$K=\Q$ and $r_2nd=n$, so that
$\frac{(n-2)(n+1)}{2}\le 1$ and one deduces that
$n=2$. Thus $A=M_2(\Q)$.

(c) Finally if $n=r_2=1$ and $d=2$ then $A$ is a
quaternion algebra over a totally real number field
which is not ramified at exactly one infinite place.
\end{proof}

The following corollary is an immediate consequence
of Propositions~\ref{ClasiKT} and \ref{vcd2}. Of
course it yields at once that (C) implies (D).

\begin{corollary}\label{vcd2KT}
Let $A$ be a finite dimensional simple rational
algebra and $R$ an order in $A$. If the virtual
cohomological dimension of $R^1$ is at most 2 then
$A$ is of Kleinian type.
\end{corollary}

\begin{remark}\label{ListOfAKT}
{\rm By Proposition~\ref{ClasiKT} there are six
types of simple algebras of Kleinian type: (1)
number fields; (2) totally definite quaternion
algebras; (3) $M_2(\Q)$; (4) $M_2(K)$, where $K$ is
an imaginary quadratic extension of the rationals;
(5) quaternion division algebras over totally real
number fields which are not ramified at exactly one
infinite place; and (6) quaternion division algebras
with exactly one pair of complex (non-real)
embeddings which are ramified at all the real
places.

Proposition~\ref{vcd2} shows that the first five
types correspond to the simple finite dimensional
rational algebras $A$ such that $\vcd(R^1)\le 2$ for
some (any) order $R$ in $A$. In the sixth case
$\vcd(R^1)=3$ (by (3.1)).

``(D) implies (E)'' of Theorem 1 (which will be
proved in sections 4, 5 and 6) shows that if $A$ is
a simple component of $\Q G$ for $G$ a finite group
of Kleinian type then $A$ is of one of the first
four types of simple Kleinian algebras.}
\end{remark}

\section{(F) implies (E)}

To prove that (F) implies (E) we need to compute the
simple components of $\Q G$, for $G$ a finite group
satisfying (F). This we will do using a method
introduced in \cite{ORS}.

Let $G$ be a finite group. For a subgroup $H$ of $G$
we set $\gorro{H}=\frac{1}{|H|}\sum_{h\in H} h$, an
idempotent element in $\Q G$. If $g\in G$ then put
$\gorro{g}=\gorro{\GEN{g}}$. A {\em strong Shoda
pair} of $G$ is a pair $(K,H)$ of subgroups of $G$
such that $H\unlhd K \unlhd G$, $K/H$ is cyclic and
$K/H$ is maximal abelian in $N_G(H)/H$. (The
definition in \cite{ORS} is more general but for our
purposes we do not need such a generality.) If $K=H$
(and hence $K=G$), then let
$\varepsilon(K,K)=\gorro{K}$; otherwise, let
$\varepsilon(K,H) = \prod_{L\in M}
(\gorro{H}-\gorro{L})$, where $M$ is the set of
minimal elements in the set of subgroups of $K$
containing $H$ properly. Finally, let $e(G,K,H)$
denote the sum of the different $G$-conjugates of
$\varepsilon(K,H)$.

Let $R$ be a ring and let $G$ be a group. If
$\rho\in \Aut(R)$ and $r\in R$ then we denote
$\rho(r)$ as $r^{\rho}$. Recall from \cite{P} that a
{\em crossed product} of $G$ over $R$ with action
$\sigma:G\rightarrow \Aut(R)$ and twisting
$\tau:G\times G \rightarrow R^*$ is an associative
ring $R*G=R*^{\sigma}_{\tau}G$ which contains $R$ as
a subring and a set of units $\{u_g \mid g\in G\}$
of $R*G$ such that $R*G=\oplus_{g\in G} u_gR$ (a
free right $R$-module) and the product in $R*G$ is
given by:
 $$(u_g r) (u_h s) = u_{gh} \tau(g,h) r^{\sigma(h)} s, \quad (g,h \in G, r,s\in R).$$

\begin{proposition}\label{e-pci-meta} \cite{ORS}
Let $G$ be a finite group.
\begin{enumerate}
\item Assume that $(K,H)$ is a strong Shoda pair of $G$. Let $N=N_{G}(H)$, $k=[K:H]$ and $n=[G:N]$.
The following properties hold.
\begin{enumerate}
\item
$e= e(G,K,H)$ is a primitive central idempotent of
$\Q G$.
\item
$\Q Ge$ is isomorphic with
$M_{n}(\Q(\xi_k)*^{\sigma}_{\tau} N/K)$, an $n\times
n$-matrix ring over a crossed product of $N/K$ over
the cyclotomic field $\Q(\xi_{k})$, with defining
action and twisting given as follows: Let $x$ be a
generator of $K/H$ and let $\gamma:N/K\rightarrow
N/H$ be a left inverse of the natural epimorphism
$N/H\rightarrow N/K$. Then
    $$\begin{array}{rl} \xi_k^{\sigma(a)} =
    \xi_k^i, & \mbox{if } x^{\gamma(a)}= x^i;\\
    \tau(a,b) = \xi_k^j, & \mbox{if }
    \gamma(ab)\inv \gamma(a)\gamma(b) = x^j,
    \end{array}$$
for $a,b\in N/K$ and integers $i$ and $j$.
\item
The simple algebra $\Q Ge$ has degree $[G:K]$.
\item
The kernel of the natural group homomorphism
$G\rightarrow Ge$ is $\core_G(H)= \bigcap_{g\in G}
H^g$.
\end{enumerate}
\item If $G$ is metabelian then every primitive central idempotent of $\Q G$ is of the form $e(G,K,H)$ for some
strong Shoda pair $(K,H)$ of $G$.
\end{enumerate}
\end{proposition}

\begin{proof}
Let $\theta$ be a linear character of $K$ with
kernel $H$. Then the induced character
$\chi=\theta^G$ is irreducible and $e=e(G,K,H)$ is
the unique primitive central idempotent of $\Q G$
such that $\chi(e)\ne 0$ \cite{ORS}. This proves
1(a). The proofs of 1(b) and 2 can be found in
\cite{ORS} and 1(c) is a direct consequence of 1(b).

To prove 1(d) note that the kernel of $g\mapsto ge$
coincides with the kernel of $\chi$. Since $H\unlhd
K\unlhd G$, this kernel is $\{k\in K :
\theta(gkg\inv) = 1, \mbox{ for all } g\in G\} =
\bigcap_{g\in G} \{k\in K : \theta(gkg\inv) = 1,
\mbox{ for all } g\in G\} = \bigcap_{g\in G} H^g =
\core_G(H)$.
\end{proof}

The following isomorphisms can be found in \cite[p.
161-163]{CR}, \cite[Lemma 20.4]{seh-book2} and
\cite{JL3}. (This can also can be  verified using
Proposition~\ref{e-pci-meta}.)

    \begin{equation}\label{WDDQ}
    \matriz{{rcl}
    \Q C_n &\cong& \oplus_{d|n} \Q(\xi_d) \\
    \Q D_{2n} &\cong& \Q (D_{2n}/D'_{2n}) \oplus \oplus_{d|n,2<d} M_2(\Q(\xi_d+\xi_d\inv)) \\
    \Q Q_{2^n} &\cong& \Q D_{2^{n-1}} \oplus \HQ(\Q(\xi_{2^{n-1}}+\xi_{2^{n-1}}\inv))\\
    \Q D_{16}^- &\cong& 4\Q \oplus M_2(\Q) \oplus M_2(\Q(\sqrt{-2})) \\
    \Q D_{16}^+ &\cong& 4\Q\oplus 2\Q(i) \oplus M_2(\Q(i))\\
    \Q \mathcal{D} &\cong& 8\Q\oplus M_2(\Q(i))\\
    \Q \mathcal{D}^+ &\cong& 4\Q\oplus 2\Q(i)\oplus 2M_2(\Q)\oplus 2M_2(\Q(i))}
    \end{equation}

Next we prove three reduction lemmas.

\begin{lemma}\label{Closed}
\begin{enumerate}
\item The class of algebras of Kleinian type is closed under
epimorphic images and semi-simple subalgebras. \item The class of
finite groups of Kleinian type is closed under subgroups and
epimorphic images. \item The class of finite groups satisfying
condition (E) of Theorem~\ref{Main} is closed under subgroups and
epimorphic images.
\end{enumerate}
\end{lemma}

\begin{proof}
1. Obviously  the class of algebras of Kleinian type is closed
under epimorphic images.

Let $A$ be an algebra of Kleinian type and $B$ a semi-simple
subalgebra of $A$. If $B_1$ is a simple quotient of $B$ then $B_1$
is isomorphic to a subalgebra of a simple quotient of $A$. In
order to prove that $B$ is an algebra of Kleinian type we thus may
assume that $A$ and $B$ are simple and $B$ is not a field. Since
$A$ is of Kleinian type and $B$ is non-abelian, it is clear that
$A$ and $B$ are quaternion algebras and that there is an order $R$
in $A$ and an embedding $\sigma:A\rightarrow M_2(\C)$ such that
$\sigma(R^1)$ is a discrete subgroup of $\SL_2(\C)$. Then $S=R\cap
B$ is an order in $B$ and $\sigma(S^1)\subseteq \sigma(R^{1})$ is
a discrete subgroup of $\SL_2(\C )$. This finishes the proof of 1.

2. This is a direct consequence of 1.

3. This can be proved imitating the proof of 1 and 2
and noticing that if $A$ and $B$ are as above ($B$
non-commutative and simple) then $Z(B)\subseteq
Z(A)$.
\end{proof}

For a finite group $G$ we denote by $\CC(G)$ the set
of isomorphism classes of noncommutative simple
quotients of $\Q G$. For simplicity we often
represent $\CC(G)$ by listing a set of
representatives of its elements. For example,
$\CC(D_{16}^+) = \{M_2(\Q(i))\}$ and $\CC(D_{16}^-)
= \{M_2(\Q),M_2(\Q(\sqrt{-2}))\}$ (see
(\ref{WDDQ})).

\begin{lemma}\label{SuperNormal}
Let $G$ be a finite group and $A$ an abelian
subgroup of $G$ such that every subgroup of $A$ is
normal in $G$. Let $\mathcal{H}=\{H\mid H \mbox{ is
a subgroup of } A \mbox{ with } A/H \mbox{ cyclic
and } G'\not\subseteq H\}$. Then
$\CC(G)=\cup_{H\in\mathcal{H}} \CC(G/H)$.
%
\end{lemma}

\begin{proof}
Let $H$ be a subgroup of $A$. By assumption,
$H\unlhd G$ and thus $\Q G = \Q G\widehat{H} \oplus
\Q G (1-\widehat{H}) \cong \Q(G/H) \oplus \Q G
(1-\widehat{H})$. It follows that $\CC(G)\supseteq
\cup_{H\in\mathcal{H}} \CC(G/H)$. It is well known
(and can be proved using
Proposition~\ref{e-pci-meta}) that the primitive
central idempotents of $\Q A$ are the elements of
the form $\varepsilon(A,H)$, where  $H$ runs through
the set $\overline{\mathcal{H}} = \{H \le A\mid A/H
\mbox{ is cyclic}\}$. Notice that each
$\varepsilon(A,H)$ is central in $\Q G$ because $H$
and $A$ are normal in $G$. Thus $\{\varepsilon(A,H)
\mid H\in \overline{\mathcal{H}}\}$ is a complete
set of orthogonal central idempotents of $\Q G$
which are primitive central in $\Q A$ but not
necessarily in $\Q G$. If $f$ is a primitive central
idempotent of $\Q G$ then there is $H\in
\overline{\mathcal{H}}$ such that
$f\varepsilon(A,H)=f$ and $f\varepsilon(A,K)=0$ for
each $H\ne K \in \overline{\mathcal{H}}$. Then $f\in
\Q G\varepsilon(A,H)$ and hence $f\in \Q G\gorro{H}
\cong \Q(G/H)$. Therefore $\Q Gf$ is a simple
epimorphic image of $\Q(G/H)$. If $\Q Gf$ is
non-commutative then $G/H$ is non-abelian and thus
$H\in\mathcal{H}$ and $\Q Gf\in \CC(G/H)$.
\end{proof}

\begin{lemma}\label{TimesAbelian}
Let $A$ be a finite abelian group of exponent $d$
and $G$ an arbitrary group.
\begin{enumerate}
\item If $d|2$ then $\CC(A\times G)=\CC(G)$.
\item If $d|4$ and $\CC(G) \subseteq \left\{M_2(\Q),\HQ(\Q), \quat{-1,-3}{\Q}, M_2(\Q(i))\right\}$ then
$\CC(A\times G)\subseteq \CC(G) \cup
\{M_2(\Q(i))\}$.
\item If $d|6$ and $\CC(G) \subseteq \left\{M_2(\Q),\HQ(\Q), \quat{-1,-3}{\Q}, M_2(\Q(\xi_3))\right\}$ then
$\CC(A\times G)\subseteq \CC(G) \cup
\{M_2(\Q(\xi_3))\}$.
\end{enumerate}
\end{lemma}

\begin{proof}
Recall that if $G_1$ and $G_2$ are two groups then
$\Q(G_1 \times G_2) \cong \Q G_1 \otimes_{\Q} \Q
G_2$. Because of the first isomorphism in
(\ref{WDDQ}), this implies, in particular, that the
simple quotients of $\Q A$ are of the form
$\Q(\xi_k)$, for $k$ a divisor of $d$. It then also
follows that the elements of $\CC(A \times H)$ are
represented by the simple quotients of the algebras
of the form $\Q(\xi_k)\otimes_{\Q} B$ for $k|d$ and
$B\in \CC(H)$. Hence, if $d=2$ each $\Q(\xi_k)=\Q$
and thus we obtain that $\CC(A\times G)=\CC(G)$. If
$d|4$ then $\Q A$ is isomorphic to a direct product
of copies of $\Q$ and $\Q(i)$ and thus if, moreover,
$\CC(G)\subseteq \left\{M_2(\Q),\HQ(\Q),
\quat{-1,-3}{\Q}, M_2(\Q(i))\right\}$ then
$\CC(A\times G)$ is formed by elements of $\CC(G)$
and simple quotients of $\Q(i)\otimes_{\Q} M_2(\Q)
\cong M_2(\Q(i))$, $\Q(i)\otimes_{\Q} \HQ(\Q)\cong
M_2(\Q(i))$, $\Q(i)\otimes_{\Q}
\quat{-1,-3}{\Q}\cong M_2(\Q(i))$ and
$\Q(i)\otimes_{\Q} M_2(\Q(i))\cong 2 M_2(\Q(i))$.
This proves 1 and 2. To prove 3 one argues similarly
using that if $d|6$ then every simple quotient of
$\Q A$ is isomorphic to either $\Q$ or $\Q(\xi_3)$
and $\Q(\xi_3)\otimes_{\Q} M_2(\Q) \cong
\Q(\xi_3)\otimes_{\Q} \HQ(\Q) \cong
\Q(\xi_3)\otimes_{\Q} \quat{-1,-3}{\Q} \cong
M_2(\Q(\xi_3))$.
\end{proof}

Now we compute $\CC(G)$ for some of the groups $G$
listed in (F) of Theorem~\ref{Main}.

\begin{lemma} \label{newlem2}
\begin{enumerate}
\item $\CC(\W_{1n}) = \{M_2(\Q)\}$.
\item $\CC(\W) = \CC(\W_{2n}) = \{M_2(\Q), \HQ(\Q)\}$.
\item $\CC(\V), \CC(\V_{1n}), \CC(\V_{2n}), \CC(\U_1), \CC(\U_2), \CC(\T_{1n}) \subseteq \{M_2(\Q), \HQ(\Q),
M_2(\Q(i))\}$.
\item $\CC(\T), \CC(\T_{2n}), \CC(\T_{3n}) \subseteq \{M_2(\Q), \HQ(\Q), M_2(\Q(i)), \HQ(\Q(\sqrt{2})),
M_2(\Q(\sqrt{-2}))\}$.
\item Let $G=M\rtimes P$ be a semidirect
product, where $M$ is a non trivial elementary
abelian $3$-group. Suppose the centralizer
$Q=\Cen_{P}(M)$ of $M$ in $P$ has index 2 in $P$ and
$m^p=m\inv$ for every $p\in P\setminus Q$.
\begin{enumerate}
\item
If $P=\GEN{x}$ is cyclic of order $2^n$ then
$\CC(G)=\CC(G/\GEN{x^2})\cup
\left\{\quat{\xi_{2^{n-1}},-3}{\Q(\xi_{2^{n-1}})}\right\}$.
In particular, if $P=C_8$ then
$\CC(G)=\left\{M_2(\Q), \quat{-1,-3}{\Q},
M_2(\Q(i))\right\}$.
\item
If $P=\W_{1n}$ and
$Q=\GEN{y_1,\ldots,y_n,t_{1},\ldots , t_{n},x^{2}}$
then $\CC(G)=\left\{M_2(\Q), \quat{-1,-3}{\Q},
M_2(\Q(\xi_{3}))\right\}$.
\item
If $P=\W_{21}$ and $Q=\GEN{y_{1}^{2},x}$ then
$\CC(G)=\{M_2(\Q), \HQ(\Q(\sqrt{3})), M_2(\Q(i)),
M_2(\Q(\xi_{3}))\}$.
\end{enumerate}
\end{enumerate}
\end{lemma}

\begin{proof}
We use the notation and presentation of the groups
as given in part (F) of Theorem~\ref{Main}. For the
groups $\W_{11}$, $\V_{11}$ and $\T_{11}$ we put
$y=y_1$ and $t=t_1$.

Throughout  the proof we will use
Lemma~\ref{SuperNormal}, Lemma~\ref{TimesAbelian}
and (\ref{WDDQ}). For several of the groups $G$
mentioned in the statement of the lemma, we will
identify a group $A$ satisfying the conditions of
Lemma~\ref{SuperNormal}. By $\mathcal{H}$ we then
denote the set (depending on $A$) considered in
Lemma~\ref{SuperNormal}. So that $\CC(G)=\cup_{H\in
\mathcal{H}} \CC(G/H)$.

\emph{1 and 2}. For $\W_{11}$, let
$A=\GEN{x^2,t}=Z(\W_{11})$. If $H\in \mathcal{H}$
then $H=\GEN{x^2}$ or $\GEN{tx^2}$ and hence
$\W_{11}/H\cong D_8$. Thus
$\CC(\W_{11})=\{M_2(\Q)\}$.

For $\W$, let $A=Z(\W)=\GEN{x^2,y^2,t}$. If $H\in
\mathcal{H}$, then $\W/H$ is a non-abelian group of
order 8 and thus it is isomorphic to $D_8$ or $Q_8$.
Hence $\CC(\W)\subseteq \{M_2(\Q), \HQ(\Q)\}$. The
converse inclusion is clear, because $D_8$ and $Q_8$
are epimorphic images of $\W$.

Since $\W_{21}$ is an epimorphic image of $\W$, and
$D_8$ and $Q_8$ are epimorphic images of $\W_{21}$,
one has that $\CC(\W_{21})=\{M_2(\Q), \HQ(\Q)\}$.

For $G=\W_{mn}$ with $m=1$ or $2$ and $n\ge 2$,
consider $A=G'=\GEN{t_1,\ldots,t_n}$. Then, using
the relations $y^2=(y,x)^{m-1}$ for every $y\in
\GEN{y_1,\dots,y_n}$, one deduces that $G/H\cong
C_2^{n-1}\times \W_{m1}$ for every $H\in
\mathcal{H}$ and thus $\CC(\W_{mn})=\CC(\W_{m1})$.

\emph{3}. For $\V$, take
$A=Z(\V)=\GEN{x^2,y^2,(y,x)}$ and let $H\in
\mathcal{H}$. If $[A:H]=2$ then $\V/H$ has order 8
and then $\CC(\V/H)$ is either $\{M_2(\Q)\}$ or
$\{\HQ(\Q)\}$. Otherwise $A/H$ is cyclic of order
$4$ and therefore $x^4\not\in H$ or $y^4\not\in H$.
Thus $\V/H$ is a group of order 16, exponent 8 and
with commutator subgroup  of order 2. This implies
that $\V/H\cong D_{16}^+$ and
$\CC(\V/H)=\{M_2(\Q(i))\}$. Thus $\CC(\V)\subseteq
\{M_2(\Q), \HQ(\Q), M_2(\Q(i))\}$.

For $G=\T_{11}$, we need a different argument.
Consider $K=\GEN{t,y,x^2}$, an abelian subgroup of
index $2$ in $G$, and the following eight subgroups
of $K$:
    $$\matriz{{llll}
    H_1 = \GEN{ y, x^2 }, &
    H_2 = \GEN{ y, x^2t^2}, &
    H_3 = \GEN{ y, tx^2 }, &
    H_4 = \GEN{ y, tx^{-2} }, \\
    H_5 = \GEN{ yt^2, x^2 }, &
    H_6 = \GEN{ yt^2, x^2t^2}, &
    H_7 = \GEN{ yt^2, tx^2 }, &
    H_8 = \GEN{ yt^2, tx^{-2} } }$$
A straightforward calculation shows that
$K=N_G(H_i)$ and $K/H_i$ is cyclic (generated by the
class of $t$) of order $4$ for every $i$, so that
$(K,H_i)$ is a strong Shoda pair of $G$ for every
$i$. By Proposition~\ref{e-pci-meta}, each
$e_i=e(G,K,H_i)=\varepsilon(K,H_i)+\varepsilon(K,H_i^x)$
is a primitive central idempotent of $\Q
G(1-\widehat{t^2})$ and $\Q Ge_i \cong M_2(\Q(i))$.
Furthermore the $16$ subgroups of the form $H_i$ and
$H_i^x$ are pairwise different. This implies that
the $8$ primitive central idempotents $e_i$ are
pairwise different and hence $\Q G=\Q G\widehat{t^2}
\oplus \oplus_{i=1}^8 \Q e_i \cong \Q(G/\GEN{t^2})
\oplus 8 M_2(\Q(i))$. Since $G/\GEN{t^2}$ is an
epimorphic image of $\V$, one concludes that
$\CC(\T_{11})\subseteq \{M_2(\Q), \HQ(\Q),
M_2(\Q(i))\}$. Actually
$\CC(\V)=\CC(\T_{11})=\{M_2(\Q), \HQ(\Q),
M_2(\Q(i))\}$ because $\W$ is an epimorphic image of
$\T_{11}$ and $\V$.

For $G=\V_{1n}, \V_{2n}, \U_1, \U_2$ or $\T_{1n}$,
we consider $A=G'$ and let $H\in \mathcal{H}$. If
$G=\T_{1n}$ then $G/H$ is an epimorphic image of
$C_4^{n-1}\times \T_{11}$, and otherwise $G/H$ is an
epimorphic image of $C_2^k\times \V$ for some $k$.
We conclude that $\CC(G)\subseteq \{M_2(\Q),
\HQ(\Q), M_2(\Q(i))\}$.

\emph{4}. For $G=\T$, take $A=Z(G)=\GEN{t^2,ty^2}$.
Let $H\in\mathcal{H}$. If either $t^2\in H$ or
$y^4\in H$ then $G/H$ is an epimorphic image of
either $\V$ or $\T_{11}$ and so $\CC(G/H)\subseteq
\{M_2(\Q), \HQ(\Q), M_2(\Q(i))\}$. Otherwise, that
is, if $t^2,y^4\not\in H$, then $t^2y^4\in H$ and
hence $H=\GEN{ty^2}$ or $H=\GEN{t\inv y^2}$. Then
$G/H$ is isomorphic to either $Q_{16}$ or $D_{16}^-$
and $\CC(G/H)\subseteq
\{M_2(\Q),\HQ(\Q(\sqrt{2})),M_2(\Q(\sqrt{-2}))\}$.
We conclude that $\CC(\T)\subseteq \{M_2(\Q),
\HQ(\Q), M_2(\Q(i)), \HQ(\Q(\sqrt{2})),
M_2(\Q(\sqrt{-2}))\}$.

For $\T_{21}$, take $A=Z(\T_{21})=\GEN{t^2,x^2}$ and
$H\in\mathcal{H}$. If $t^2\in H$ or $t^2x^2\in H$
then $G/H$ is an epimorphic image of $\V$ or $\T$.
Otherwise, $H=\GEN{x^2}$ and hence $G/H=D_{16}^-$.
So $\CC(\T_{21})\subseteq \{M_2(\Q), \HQ(\Q),
M_2(\Q(i)), \HQ(\Q(\sqrt{2})),
M_2(\Q(\sqrt{-2}))\}$.

For $G=\T_{2n}$, taking $A=G'$ and having in mind
that $(y,x)y^2 = 1$ for every $y\in
\GEN{y_1,\ldots,y_n}$, one deduces that $G/H$ is an
epimorphic image of $T_{21}\times C_2^{n-1}$ for
every $H\in \mathcal{H}$. Hence, by the previous
paragraph, $\CC(\T_{2n})=\CC(\T_{21})$.

Finally, for $G=\T_{3n}$, take $A=G'$ and let $H\in
\mathcal{H}$. If $t_1^2\in H$ then $G/H$ is an
epimorphic image of $\T_{3n}/\GEN{t_1^2}$ which in
turn is an isomorphic image of $\V_{1n}$ and thus
$\CC(G/H)\subseteq \{M_2(\Q), \HQ(\Q),
M_2(\Q(i))\}$. Otherwise, the image $\overline{t_i}$
in $G/H$ of each $t_i$ belongs to
$\GEN{\overline{t_1}}$. Furthermore,
$\overline{t_1}$ has order 4 and $\overline{t_i}$
has order at most 2, for $i\ge 2$. Thus
$\overline{t_i}\in \GEN{t_1^2}$. For $i\ge 2$, let
$z_i$ be the natural image of $y_i$ in $G/H$ if
$t_i=1$ and the natural image of $y_1^2y_i$,
otherwise. Then $Z=\GEN{z_2,\ldots,z_n}$ is central
in $G/H$ of exponent at most 2 and $G/H$ is an
epimorphic image of $\T\times Z$. Thus
$\CC(G/H)\subseteq \CC(\T)$.

\emph{5}. Let $G=M\rtimes P$ be as in statement 5.
Applying Lemma~\ref{SuperNormal}, with $A=M$, we may
assume, without loss of generality, that $M$ is
cyclic of order 3, generated by $m$, say.

(a) Assume $P=\GEN{x}$ is cyclic of order $2^n$.
Then $(K=\GEN{m,x^2},1)$ is a strong Shoda pair of
$G$ and so $e=e(G,K,1)$ is a primitive central
idempotent of $\Q G$. Applying
Proposition~\ref{e-pci-meta} one has $\Q
Ge=\Q(\xi)[u|u^2=\xi^3,u\inv \xi u = \xi^s]$, where
$\xi=\xi_{3\cdot 2^{n-1}}$ and $s$ is an integer
such that $s\equiv -1 \mod 3$ and $s\equiv 1 \mod
2^{n-1}$. Let $\omega=\xi^{2^{n-1}}$, a third root
of unity. Then $j^2=-3$ and $ju=-uj$, where
$j=1+2\omega$. This shows that $\Q Ge\cong
\quat{\xi_{2^{n-1}},-3}{\Q(\xi_{2^{n-1}})}$. Since
$e+(1-\widehat{G'})\widehat{x^2}=1-\widehat{G'}$,
one concludes that $\CC(G)=\CC(G/\GEN{x^2})\cup\{\Q
Ge\}$.

In particular, if $n=1$ then
$\CC(G)=\{\quat{1,-3}{\Q}=M_2(\Q)\}$, if $n=2$ then
$\CC(G)=\{M_2(\Q),\quat{-1,-3}{\Q}\}$ and if $n=3$
then $\CC(G)=\{M_2(\Q),\quat{-1,-3}{\Q},
\quat{i,-3}{\Q(i)}=M_2(\Q(i))\}$. The equality
$M_{2}(\Q (i)) = \quat{i,-3}{\Q(i)}$ holds because
the equation $iX^2-3Y^2=1$ has the solution $X=1+i$
and $Y=i$.

(b) Assume $P=\W_{1n}$. Applying
Lemma~\ref{SuperNormal} with $A=P'$, we may assume
that $n=1$, because $G/H\cong C_2^{n-1}\times
(M\rtimes \W_{11})$ for every $H\in \mathcal{H}$. So
suppose $P=\W_{11}$ and $Q=\GEN{x^2,y=y_1,t=t_1}$.
Let $S$ be a simple quotient of $\Q G$. Put
$T=G/\GEN{t}$ and
$e=(1-\widehat{m})(1-\widehat{t})$, a central
idempotent of $\Q G$. Notice that $G/\GEN{m}\cong P$
and $T\cong C_2\times (C_3\rtimes C_4)$. If $S$ is a
quotient of $\Q G(1-e)$ then $S$ is a quotient of
either $\Q G \widehat{m} \cong \Q(G/\GEN{m})\cong \Q
P$ or $\Q G \widehat{t} \cong \Q T \cong
\Q(C_2\times (C_3\rtimes C_4))$. Then $S\cong
M_2(\Q)$ or $\quat{-1,-3}{\Q}$, by the previous
paragraph and statement 1. Otherwise, $S$ is a
quotient of $\Q Ge$. Notice that
$e=e(G,K,H_1)+e(G,K,H_2)$, where $K=\GEN{Q,m}$,
$H_1=\GEN{x^2,y}$ and $H_2=\GEN{tx^2,y}$ and
$(K,H_1)$ and $(K,H_2)$ are two strong Shoda pairs
of $G$. Clearly $[K:H_1]=[K:H_2]=6$ and
$K=N_G(H_1)=N_G(H_2)$. Hence, because of
Proposition~\ref{e-pci-meta}, it follows that $\Q Ge
\cong 2M_2(\Q(\sqrt{-3}))$, and therefore $S\cong
M_2(\Q(\sqrt{-3}))$.

(c) Assume $P=\W_{12}$ and $Q=\GEN{x,y^2}$. Let $S$
be a simple quotient of $\Q G$. Again, put
$T=G/\GEN{t}$ and
$e=(1-\widehat{m})(1-\widehat{t})$. Then
$G/\GEN{m}\cong P$ and $T\cong C_4\times (C_3\rtimes
C_2)$. Thus, if $S$ is a quotient of $\Q G(1-e)$
then $S$ is a quotient of either $\Q G \widehat{m}
\cong \Q(G/\GEN{m})\cong \Q P$ or $\Q G \widehat{t}
\cong \Q T \cong \Q(C_4\times (C_3\rtimes C_2))$.
Then, as in the proof of (b), $S$ is isomorphic to
either $M_2(\Q)$, $\HQ(\Q)$ or $M_2(\Q(i))$.
Otherwise, $S$ is a quotient of $\Q Ge$. In this
case, $e=e(G,K,H_1)+e(G,K,H_2)$, where
$K=\GEN{Q,m}$, $H_1=\GEN{x}$, $H_2=\GEN{x^2y^2}$,
and $(K,H_1)$ and $(K,H_2)$ are two strong Shoda
pairs of $G$. Since $[K:H_1]=6$ and $K=N_G(H_1)$ we
get that  $\Q Ge(G,K,H_1)\cong M_2(\Q(\sqrt{-3}))$,
as desired. Because $N_G(H_2)=G$, we deduce from
Proposition~\ref{e-pci-meta} that $\Q Ge(G,K,H_2)$
is the simple algebra given by the following
presentation: $S=\Q(\xi)[u|u^2=-1, u\inv \xi u =
\xi\inv]$, with $\xi = \xi_{12}$. Let $i=\xi^3$,
$j=u$ and $a=\xi+\xi\inv\in Z(S)$. Clearly
$i^2=u^2=-1$, $a^2=3$ and $ji=-ij$. Therefore
$S\cong \HQ(\Q(\sqrt{3}))$. Thus $S\cong
M_2(\Q(\sqrt{-3}))$ or $\HQ(\Q(\sqrt{3}))$.
\end{proof}

We are ready to prove (F) implies (E). So, let $G$
be a finite group satisfying (F). By
Lemma~\ref{Closed}, to prove that $G$ satisfies (E)
one may assume that $G=A\times H$ for $A$ and $H$
satisfying one of the conditions 1-4 in (F). We have
to show that the elements of $\CC(G)$ either are
totally definite quaternion algebras or are of the
form $M_2(\Q(\sqrt{-d}))$ for $d=0,1,2$ or $3$.
Using Lemma~\ref{TimesAbelian}  and
Lemma~\ref{newlem2}, one obtains the following five
statements, and hence the result follows. If either
condition 1 or condition 2 holds then
$\CC(G)\subseteq \{M_2(\Q), \HQ(\Q), M_2(\Q(i))\}$.
If condition 3 holds then $\CC(G) \subseteq
\{M_2(\Q), \HQ(\Q), M_2(\Q(i)), \HQ(\Q(\sqrt{2})),
M_2(\Q(\sqrt{-2}))\}$. If condition 4 holds then
$\CC(G)$ is contained in either $\left\{M_2(\Q),
\quat{-1,-3}{\Q}, M_2(\Q(i))\right\}$,
$\left\{M_2(\Q), \quat{-1,-3}{\Q},
M_2(\Q(\sqrt{-3}))\right\}$ or $\{M_2(\Q),
\HQ(\Q(\sqrt{3})), M_2(Q(i)), M_2(\Q(\sqrt{-3}))\}$,
depending on the respective cases.

\section{(D) implies (F), for nilpotent groups} \label{nilpotent-case}

In the remainder of the paper we prove (D) implies
(F), or equivalently we classify the groups of
Kleinian type as the epimorphic images of the groups
listed in (F). In this section we do this for finite
groups that are nilpotent.

We start with two lemmas which provide information
on the groups of Kleinian type.

\begin{lemma} \label{SSP-index}
Let $G$ be a finite non-abelian group of Kleinian
type. The following properties hold.
\begin{enumerate}
\item\label{TAmitsur} Either $G/Z(G)$ is elementary abelian of order 8 or $G$ has an abelian normal subgroup of index $2$.
    In particular, $G$ is metabelian and has a nilpotent subgroup of index at most $2$.
\item\label{Indices} Every primitive central
idempotent of $\Q G$ is of the form $e = e(G,K,H)$
for some strong Shoda pair $(K,H)$ of $G$. Moreover,
for such a primitive central idempotent $e$ one has
    \begin{enumerate}
    \item $[G:K]\le 2$;
    \item if $H$ is not normal in $G$ then
    $K=N_G(H)$ and $[K:H]$ divides $4$ or
    $6$, and
    \item if $\Q G e$ is not a division ring then $[K:H]$ divides $8$ or $12$ and
    $\Q Ge$ is isomorphic to $M_2(\Q(\sqrt{-d}))$ for $d=0,1,2$ or $3$.
    \end{enumerate}
\item\label{Dihedral}
If a dihedral group $D_{2n}$ is an epimorphic image
of a subgroup of $G$ then $n$ divides $4$ or $6$.
\item\label{Descomposicion}
$G=G_3\rtimes G_2$ where $G_3$ is an elementary
abelian $3$-group (possibly trivial), $G_2$ is a
$2$-group and the kernel of the action of $G_2$ on
$G_3$ has index at most 2 in $G_2$.
\item\label{ExpCentro} The exponent of $Z(G)$ is a divisor of $4$ or $6$.
\item\label{Derivado}
$Z(G)\cap G'=\{t\in G'\mid t^2=1\}$. Furthermore, if
$t\in G'$ and $x\in G$ then either $t^x=t$ or
$t^x=t\inv$. In particular, every subgroup of $G'$
is normal in $G$.
\end{enumerate}
\end{lemma}

\begin{proof}
1. Since every simple quotient of $\Q G$ has degree
at most $2$ (see the definition of groups and
algebras of Kleinian type), the irreducible
character degrees of $G$ are $1$ and $2$. By
\cite{AMI} this implies that either $G/Z(G)$ is
elementary abelian of order $8$ or $G$ has an
abelian subgroup of index $2$. In the first case $G$
is central-by-abelian, and hence nilpotent and
metabelian. In the second case, obviously $G$ is
also metabelian and it has a nilpotent (in fact
abelian) subgroup of index $2$.

2. That every primitive central idempotent of $\Q G$
is of the form $e=e(G,K,H)$ for some strong Shoda
pair of $G$ is a consequence of
Proposition~\ref{e-pci-meta} and the fact that $G$
is metabelian. Let $e=e(G,K,H)$ for $(K,H)$ a strong
Shoda pair of $G$.

The inequality $[G:K]\le 2$ is a straightforward
consequence of the fact that $[G:K]$ equals the
degree of $\Q Ge$ (Proposition 4.1). Let $k=[K:H]$.
Since $K\le N_G(H)$, we get that either $K=N_G(H)$
or $G=N_G(H)$. Therefore, if $H$ is not normal in
$G$ then $K=N_G(H)$ and $\Q Ge=M_2(\Q(\xi_k))$. By
Proposition~\ref{ClasiKT},
$\varphi(k)=[\Q(\xi_k):\Q]\le 2$ and therefore $k$
divides $4$ or $6$. This proves (b) and it also
proves (c) if $H$ is not normal in $G$. If $\Q Ge$
is not a division ring and $H$ is normal in $G$
then, by Proposition~\ref{e-pci-meta}, $\Q Ge\cong
M_2(F)$ where $F$ is a subfield of index 2 in
$\Q(\xi_k)$. Furthermore, Remark~\ref{ListOfAKT}
implies that $F$ is either $\Q$ or an imaginary
quadratic extension of $\Q$. Hence
$\varphi(k)=2[F:\Q]$, a divisor of 4. If
$\varphi(k)\ne 4$ then $k$ divides $4$ or $6$ and
therefore $A=M_2(\Q)$. Otherwise $k=5,8,10$ or $12$.
If $5|k$ then necessarily $F=\Q(\sqrt{5})$, a
contradiction. Thus $k=8$ or $12$ and therefore
$A=M_2(\Q(\sqrt{-d}))$ for $d=1,2$ or $3$. This
finishes the proof of 2.

3. By (\ref{WDDQ}), $\Q D_{2n}$ has an epimorphic
image isomorphic to $M_2(\Q(\xi_n+\xi_n\inv))$. This
algebra is of Kleinian type, by Lemma~\ref{Closed}.
Therefore $\Q(\xi_n+\xi_n\inv)=\Q$, by
Remark~\ref{ListOfAKT} and this implies that $n$
divides $4$ or $6$.

4. Let $E$ be the set of primitive central
idempotents $e$ of $\Q G$ such that $\Q Ge$ is
non-commutative. Let $e\in E$ and $z\in Z(G)$. By
Proposition~\ref{ClasiKT} the order of $ze$ divides
$4$ or $6$ and thus the exponent of $Z(G)$ divides
$12$.

By 1, $G$ has a nilpotent subgroup of index at most
2. Therefore, $G=G_{2'} \rtimes G_2$, where $G_{2'}$
is a nilpotent subgroup of odd order of $G$, $G_2$
is a Sylow $2$-subgroup of $G$ and the kernel of the
action of $G_2$ on $G_{2'}$ has index at most $2$ in
$G_2$. We have to show that $G_{2'}$ is an
elementary abelian $3$-group. We argue by
contradiction. So, let $a\in G_{2'}$ be of order
$q$, where $q$ is either $9$ or $q>3$ and prime.
Since the exponent of $Z(G)$ divides $12$, $a$ is
not central in $G$ and so there is $x\in G_2$ such
that $a^x\ne a$. Put $b=aa^x$. Assume that $b=1$.
Then $\GEN{a,x}/\GEN{x^2} \cong D_{2q}$ and thus
$D_{2q}$ is of Kleinian type by Lemma~\ref{Closed},
contradicting \ref{Dihedral}. Therefore $b$ is a
non-trivial central element of odd order. Hence $b$
has order $3$ and $a$ has order $9$. If $b\in
\GEN{a}$, then $b=a^{\pm 3}$ and hence $a^x=a^2$ or
$a^x=a^{-4}$. Then $a=a^{x^2} = a^4$ or
$a=a^{x^2}=a^7$, a contradiction. Thus
$\GEN{a,x}/\GEN{x^2}=(\GEN{a}_9\times
\GEN{b}_3)\rtimes \GEN{x}_2$, with $a^x = a\inv b$
and $b^x=b$. Then $\GEN{a,x}/\GEN{b,x^2} \cong
D_{18}$, again a contradiction.

5. Since the exponent of $Z(G)$ divides $12$ it is
enough to show that $Z(G)$ does not have elements of
order $12$. By means of contradiction assume that
$a\in Z(G)$ has order $12$. Let
$\varepsilon=\varepsilon(\GEN{a},1)$ (see the
notation introduced in Section 4). If
$\varepsilon(1-\widehat{G'})\ne 0$ then there is a
(necessarily injective) non-zero homomorphism
$\Q(\xi_{12}) \cong \Q\GEN{a} \varepsilon
\rightarrow Z(A)$ for some non-commutative simple
quotient $A$ of $\Q G$. This implies that $A$ has a
central root of unity of order $12$, contradicting
Proposition~\ref{e-pci-meta}. Thus $\varepsilon =
\varepsilon \widehat{G'}$. If $H=\GEN{a} \cap G'$
then $0\ne \varepsilon = \varepsilon \widehat{G'} =
\varepsilon \widehat{H} \widehat{G'}$. If $H\ne 1$
then $0\ne \varepsilon
\widehat{H}=(1-\widehat{a^4})(1-\widehat{a^6})=0$
because $H$ contains either $a^4$ or $a^6$. Thus
$H=1$ and so $\varepsilon  = \varepsilon
\widehat{G'} = \frac{1}{|G'|} \varepsilon \sum_{g\in
G'} g$. Since all the elements of $G$ with non-zero
coefficient in $\varepsilon$ belongs to $\GEN{a}$,
the last formula implies that $G'\subseteq \GEN{a}$.
Thus $G'=1$, contradicting the fact that $G$ is
non-abelian.

6. If $G/Z(G)$ is elementary abelian then
$G'\subseteq Z(G)$. It follows that, for $a,b\in G$
we get that $b^2a=ab^2 = batb = b^2at^2$, where
$t=(a,b)$. Then $t^2=1$ and the statement follows.

So assume that $G/Z(G)$ is not elementary abelian.
From 1, we then have that $G$ has an abelian
subgroup $N$ of index 2. Let $a\in G\setminus N$.
Then $t^x=t$ if $x\in N$. If $x\in G\setminus N$
then $x=na$ for some $n\in N$ and therefore
$t^x=t^a$. Moreover $a^2\in Z(G)$ and then, for
every $g\in G$, one has
$1=(g,a^2)=(g,a)(g,a)^{a\inv}=(g,a)(g,a)^a$. Thus
$(g,a)^a=(g,a)\inv$. On the other hand if $n,m\in N$
then $(na,ma)=(na,m)(na,a)^{m\inv} =
(a,m)^{n\inv}(n,m)\left((a,a)^{n\inv} (n,a)
\right)^{m\inv} = (a,m)(n,a)=(m,a)\inv (n,a)$. Thus
if $t\in G'$ then $t=(n_1,a)^{\alpha_1}\cdots
(n_k,a)^{\alpha_k}$ for some $n_i\in N$ and
$\alpha_i\in \Z$ and $t^x=(n_1,a)^{-\alpha_1}\cdots
(n_k,a)^{-\alpha_k}=t\inv$. So we have shown that
$t^x=t$ if $x\in N$ and $t^x = t\inv$ otherwise.
Therefore $t\in Z(G)$ if and only if $t^2=1$.
\end{proof}

\begin{lemma}\label{PropEl}
If $G$ is a finite non-abelian $2$-group of Kleinian
type then the following properties hold.
\begin{enumerate}
\item\label{Exponente} The exponent
of $G$ is at most $8$.
\item\label{ExpDerivado}
$G'$ is abelian of exponent at most $4$.
\item\label{DeriDeri} $\GEN{(G,G')}=
{G'}^2\subseteq Z(G)$.
\item\label{Derivado1} $G_x = \GEN{(x,g) \mid g\in G}$ is a normal subgroup of $G$ for all $x\in G$. Moreover, if
$G'\neq G_x $ then $x^{4}\in G_{x}$.
\item\label{DosCasos} If
$x,y\in G$ and $t=(y,x)$ then one of the following
conditions holds:
\begin{enumerate}
\item $(x,t)=1$, $(y,x^2)=t^2$ and $(y,tx^2)=1$ or
\item $(x,t)\ne 1$ and $(y,x^2)=1$.
\end{enumerate}
\end{enumerate}
\end{lemma}

\begin{proof}
Let $E$ be the set of primitive central idempotents
$e$ of $\Q G$ such that $\Q Ge$ is non-commutative.
It is well known that $1-\widehat{G'}=\sum_{e\in E}
e$ (see for example \cite{Coleman}). Notice that the
coefficient of $1$ in $g(1-\widehat{G'})$ is $0$, if
$g\not\in G'$, and it is $-\frac{1}{|G'|}$, if $1\ne
g\in G'$. However the coefficient of $1$ in
$1-\widehat{G'}$ is $1-\frac{1}{|G'|}>0$ because, by
assumption, $G'$ is non trivial. This shows that the
natural group homomorphism $G \rightarrow
G(1-\widehat{G'})$, mapping $g\in G$ onto
$g(1-\widehat{G'})$, is injective, and hence so is
the natural group homomorphism $f:G\rightarrow
\prod_{e\in E} Ge$.

1. We prove the statement by contradiction. So
suppose $G$ is a non-abelian $2$-group of Kleinian
type of minimal order such that the exponent of $G$
is greater than $8$. Let $g\in G$ be of order $16$.

Then, there is $e\in E$ such that $ge$ has order
$16$ and, by the minimality of $G$, $G$ is
isomorphic to $Ge$. By Proposition~\ref{e-pci-meta},
there is a strong Shoda pair $(K,H)$ of $G$ such
that $e=e(G,K,H)$, $[G:K]=2$ and $\core_G(H)=1$.
Then $A=\Q Ge$ has a subfield isomorphic to
$\Q(\xi_{16})$. Since $A$ is a quaternion algebra,
the dimension of the centre of $A$ is at least
$\varphi(16)/2=4$. Hence, by statement~\ref{Indices}
of Proposition~\ref{SSP-index}, $A$ is a division
algebra. It follows from statement (b) of
Theorem~\ref{e-pci-meta} implies that $H\unlhd G$,
that is, $H=\core_G(H)=1$. Thus $K$ is a cyclic
subgroup of index $2$ in $G$.
%
%
Hence, as mentioned in the preliminaries, $G$ is
isomorphic to either $D_{2^{k+1}}$, $D_{2^{k+1}}^+$,
$D_{2^{k+1}}^-$, or $Q_{2^{k+1}}$. Since $A$ is a
non-commutative division ring containing
$\Q(\xi_{16})$, $G=Q_{2^{k+1}}$, one has that $k\ge
4$ (see (\ref{WDDQ})). Thus $D_{16}$ is a quotient
of $G$, in contradiction with
statement~\ref{Dihedral} of Lemma~\ref{SSP-index}.

2. That $G'$ is abelian is a consequence of
statement~\ref{TAmitsur} of Lemma~\ref{SSP-index}.
We prove by contradiction that $G'$ has exponent at
most $4$. So, assume that $G$ is a non-abelian
$2$-group of Kleinian type of minimal order with a
commutator $t=(y,x)$ of order greater than 4. By the
minimality of the order of $G$, one has that
$G=\GEN{x,y}$. By statement \ref{Derivado} of
Lemma~\ref{SSP-index}, $t\in G'\setminus Z(G)$.
Hence $G/Z(G)$ is non-abelian and
statement~\ref{TAmitsur} of Lemma~\ref{SSP-index}
implies that $G$ has an abelian subgroup $A$ of
index 2. Therefore, either $x\not\in A$ or $y\not\in
A$. Since $(yx,x)=t$, one may assume that $x\not\in
A$ and $y\in A$. Then $xy\not\in A$ and therefore
$x^2,(xy)^2\in Z(G)$. Furthermore, by statement
\ref{Derivado} of Lemma~\ref{SSP-index}, $t^x =
t\inv$. Hence $(xy)^2=t\inv x^2 y^2$ and thus $t\inv
y^2\in Z(G)$. So, by statement~\ref{ExpCentro} of
Lemma~\ref{SSP-index}, $t^{-4}y^8=1$. Thus
$y^8=t^4\ne 1$, in a contradiction with 1.

3. $G'^2\subseteq Z(G)$ is a consequence of 2 and
statement~\ref{Derivado} of Lemma~\ref{SSP-index}.
Furthermore for $t\in G'$,  either $t^2=1$, or $t$
has order $4$ and $t\not\in Z(G)$. Thus, again by
statement~\ref{Derivado} of Lemma~\ref{SSP-index},
there is $x\in G$ such that $t^x=t\inv$, that is,
$(x,t)=t^2$. Hence 3 follows.

4. That $G_x$ is normal in $G$ is a direct
consequence of statement~\ref{Derivado} in
Lemma~\ref{SSP-index}. Clearly, the natural image of
$x$ in $G/G_x$ is central. Since $G/G_x$ is a
$2$-group of Kleinian type, it follows from
statement~\ref{ExpCentro} of Lemma~\ref{SSP-index}
that if $G_{x}$ is properly contained in $G'$ then
$x^{4}\in G_{x}$ as desired.


5. Let $x,y\in G$ and $t=(y,x)$. Clearly $(y,x^2) =
(y,x)(y,x)^{x\inv}=tt^{x\inv}$. Because of
statement~\ref{Derivado} in Lemma~\ref{SSP-index},
we also know that $t^{x\inv}=t$ or
$t^{x\inv}=t^{-1}$. In the latter case we get that
$(y,x^{2})=1$ and so (b) holds. In the former case
$(t,x)=1$ and $(y,x^2)=t^2$. If also $t^{y}=t$ then
$t$ is central in $G$, and thus, again by
statement~\ref{Derivado} in Lemma~\ref{SSP-index},
$t^{2}=1$ and $(y,t)=t^{2}=1$. If, on the other
hand, $t^{y}\neq t$, then, again by
statement~\ref{Derivado} in Lemma~\ref{SSP-index},
$(y,t)=t^{2}$.  So, in all cases we get that
$(y,x^{2})=t^{2}$. By part \ref{ExpDerivado}, we
also know that $t^{4}=1$. Hence,
$(y,tx^{2})=(y,t)(y,x^{2})^{t\inv}=(y,t)t^{2}=1$, as
desired.
\end{proof}

The next three lemmas contain information on two and
three generated 2-groups with a commutator of order
4.

\begin{lemma}\label{TodosLosCasosRango2}
Let $G=\GEN{x,y}$ be a non-abelian 2-group and
suppose $t=(y,x)$ has order $4$. If $G$ is of
Kleinian type then one of the following conditions
holds.
\begin{enumerate}
\item
$(x,t)=(y,t)=t^2$, $(xy,t)=1$,
$Z(G)=\GEN{t^2,x^2,y^2}$, and one of the following
conditions holds:
\begin{enumerate}
\item $t^2=x^4y^4$;
\item $t^2\in \{x^2,y^2,x^4y^2,x^2y^4\}$;
\item $x^4=y^4=1$ and $t^2=x^2y^2$.
\end{enumerate}
\item
$(x,t)=t^2$, $(y,t)=1$,
$Z(G)=\GEN{t^2,x^2,(xy)^2}=\GEN{t^2,x^2,ty^2}$, and
one of the following conditions holds:
\begin{enumerate}
\item $y^4=1$;
\item either $ty^2\in \{x^2,x^{-2}\}$ or $x^2\in
\{t^2,y^4\}$;
\item $x^4=1$ and $t=y^{-2}$.
\end{enumerate}
\item[2'] $(x,t)=1$, $(y,t)=t^2$, $Z(G)=\GEN{t^2,y^2,(xy)^2}=\GEN{t^2,x^2,ty^2}$, and one
of the following conditions holds:
\begin{enumerate}
\item $x^4=1$;
\item either $t\inv x^2\in \{y^2,y^{-2}\}$ or $y^2\in
\{t^2,x^4\}$;
\item $y^4=1$ and $t=x^{2}$.
\end{enumerate}
\end{enumerate}
\end{lemma}

\begin{proof} Since by assumption, $t$ has order
$4$, statement~\ref{Derivado} of
Lemma~\ref{SSP-index} yields that $t$ is not
central, $t^2\in Z(G)$ and $\GEN{t}$ is normal in
$G$.  Furthermore, $(x,t)=t^2$ or $(y,t)=t^2$. We
deal with three mutually exclusive cases.

(1) First assume $(x,t)=(y,t)=t^{2}$. So  $x t =
t\inv x$ and $y t = t\inv y$.  Hence $(y,x^2)=
(y,x)(y,x)^{x^{-1}}=tt^{x^{-1}}= 1$.  Similarly
$(x,y^2)=1$. Therefore $\langle
t^{2},x^2,y^2\rangle \subseteq Z(G)$. Since
$t\not\in Z(G)$, and thus $t\not\in \langle
t^{2},x^{2},y^{2}\rangle$, it follows that
$G/\langle t^{2},x^{2},y^{2}\rangle$ is a
non-abelian group of order $8$. Hence we obtain
that $Z(G)=\GEN{t^2,x^2,y^2}$. Also $(xy,t)=1$
because $(xy,t)=(y,t)^{x\inv}
(x,t)=(t^{2})^{x\inv}t^{2}=1$.

Let $H=\GEN{x^2,y^2}$, a central subgroup of $G$.
Note that $(xy)^{2}=t^{3}x^{2}y^{2}$ and thus, in
the group $G/H$,  one has that $\overline{t}^{-1}
=\left( \overline{xy}\right)^{2}$ and
$\left(\overline{xy}\right)^{\overline{y}}=
\left(\overline{xy}\right)^{-1}$. By
statement~\ref{Exponente} of Lemma~\ref{PropEl} we
know that $(xy)^{8}=1$. Hence, there is an
epimorphism $f:D_{16} \rightarrow G/H$  given by
$f(a)=\overline{xy}$ and $f(b)=\overline{y}$, where
$D_{16}=\GEN{a}_8\rtimes \GEN{b}_2$ is the dihedral
group of order $16$. Because of
statement~\ref{Dihedral} in Lemma~\ref{SSP-index},
$D_{16}$ is not of Kleinian type. However, by
Lemma~\ref{Closed}, $G/H$ is of Kleinian  type.
Hence, $\ker f \ne 1$ and therefore the order of
$G/H$ divides $8$. Thus
\begin{equation}\label{AlaDos} 1\ne t^2 \in
\GEN{x^2,y^2}.
\end{equation}
We consider three cases: (i) $x^{4}\neq 1$, (ii)
$y^{4}\neq 1$ and (iii) $x^{4}=y^{4}=1$.

(i) Suppose $x^{4}\ne 1$. So, by statement~\ref{Exponente} of
Lemma~\ref{PropEl}, $x^{4}$ has order $2$ and, by the above,
$x^{2}$ is central in $G$. Then $e= \frac{1}{2} (1-t^{2})\;
\frac{1}{2} (1-x^{4})$ is a nonzero central idempotent of $\Q G$.
Clearly, the semi-simple $\Q$-algebra $A=\Q Ge$ is  contained in
$\Q G (1-\widehat{t})$. The latter, and thus also $A$, is a direct
sum of non-commutative simple algebras (see the beginning of the
proof of Lemma~\ref{PropEl}).

Let $f=\widehat{x^{2}y^{2}}e$, a central idempotent
of $A$. $f=0$. Because $\{ xy , txy\}$ is a full
conjugacy class of $G$, we get that $z=(1+t)xyf$
and $i=x^{2}f$ are central elements of $Af$. Since
$x^{4}e=-e$, $t^{2}e=-e$, $x^{2}y^{2}f=f$ and
$ef=f$, we get that $i^{2}=-f$ and
$z^{2}=(1+t)^{2}(xy)^{2}f
=(1+2t+t^{2})t^{-1}x^{2}y^{2}=2f$. If $f\neq 0$
then there exists a primitive central idempotent
$f_1$ of $Af$ such that $i^2f_1=-f_1$ and
$z^2f_1=2f_1$. Then  $\Q G f_1$ is a
non-commutative simple quotient of $\Q G$ having a
central subfield isomorphic to
$\Q(\sqrt{2},i)=\Q(\xi_8)$, contradicting the last
statement of Proposition~\ref{ClasiKT}. Thus $f=0$
and this implies that $\GEN{x^{2}y^{2}}\cup
t^{2}x^{4} \GEN{x^{2}y^{2}} =
t^{2}\GEN{x^{2}y^{2}}\cup x^4 \GEN{x^{2}y^{2}}$.
Therefore
\begin{equation} \label{AlaDos1}
t^{2} \in \GEN{x^{2}y^{2}} \mbox{ or } x^{4} \in
\GEN{x^{2}y^{2}} .
\end{equation}
If $x^{4}=y^{4}$ then $(x^{2}y^{2})^{2}=x^{8}=1$,
 and it follows that $t^{2}=x^{2}y^{2}$ or
$x^{2}=y^{2}$. In both situations the central
elements $i=x^2e$ and $z=(1+t)xye$  of $A$ are such
that  $i^{2}=-e$ and $z^{2}=-2e$. Hence $A$ has a
central subfield isomorphic with $\Q (i,
\sqrt{-2})=\Q(\xi_8)$, again yielding a
contradiction. Thus $x^{4}\neq y^{4}$. Since, by
statement~\ref{Exponente} of Lemma~\ref{PropEl},
$x^{8}=y^{8}=1$, we obtain that $x^{2}y^{2}$ has
order $4$. As also both $t$ and $x^{2}$ have order
$4$, (\ref{AlaDos1}) then implies that either
$t^{2}=x^{4}y^{4}$ or $y^{4}=1$. In the first case
1(a) holds. In the second case $y^{4}=1$ and thus
(\ref{AlaDos}) implies that
$t^{2}\in\{x^4,y^2,x^4y^2\}$. If $t^{2}=x^{4}$ then
1(a) holds, otherwise 1(b) holds.

(ii) Suppose $y^{4}\ne 1$. By symmetry with case
(i), we also obtain that either 1(a) or 1(b) holds.

(iii) Suppose $x^{4}=y^{4}=1$. By (\ref{AlaDos}),
we clearly get that $x^{2}\neq 1$ or $y^{2}\neq 1$.
If $y^{2}\in \GEN{x^{2}}$ then, by (\ref{AlaDos}),
we obtain that $1\neq t^{2}\in \GEN{x^{2}}$ and
thus $t^{2}=x^{2}$; hence 1(b) holds. Similarly, if
$x^{2}\in \GEN{y^{2}}$ then $t^2=y^2$ and again
1(b) holds. Otherwise $\GEN{x^{2},y^{2}}\cong
C_{2}\times C_{2}$ and thus by (\ref{AlaDos}) one
of the following holds: $t^{2}=y^{2}$,
$t^{2}=x^{2}$ or $t^{2}=x^{2}y^{2}$. Hence 1(b) or
1(c) holds. This finishes the proof when
$(x,t)=(y,t)=t^{2}$.

(2) Second assume that $(y,t)=1$.  Hence, since $t$
is not central, $(x,t)=t^2$. Set $x_{1}=x$,
$y_{1}=yx$ and $t_{1}=(y_{1},x_{1})$. Then
$t_{1}=t$ and
$(x_{1},t_{1})=(y_{1},t_{1})=t_{1}^{2}$. Thus, by
(1), $x_{1}$ and $y_{1}$ satisfy one of the
conditions of (1).  In particular
$Z(G)=\GEN{t^2,x_1^2,y_1^2}$. Notice that
$y_{1}^{2} = t\inv x^{2} y^{2}$ and $y_{1}^{4} =
t^{2}x^{4}y^{4}$. Then $Z(G)=\GEN{t^2,x_1^2,y_1^2}
= \GEN{t^2,x^2,(xy)^2} = \GEN{t^2,x^2,ty^2}$.

Thus if $x_{1}$ and $y_{1}$ satisfy 1(a) then
$y^4=1$, that is, condition 2(a) holds.

Next assume that $x_{1}$ and $y_{1}$ satisfy 1(b).
If $t_{1}^{2}=x_{1}^{2}$ then $t^{2}=x^{2}$. If
$t_{1}^{2}=y_{1}^2$ then $t=x^{-2}y^{-2}$. {If
$t_{1}^{2}=x_{1}^{4}y_{1}^{2}$ then}
$t=x^{2}y^{-2}$. If $t_{1}^{2}=x^{2}y_{1}^{4}$ then
$x^{2}=y^{4}$. So always 2(b) holds.

Finally assume that $x_{1}$ and $y_{1}$ satisfy
1(c). Then $x^{4}=1$ and $t=y^{-2}$, that is 2(c)
holds.

(3) Third assume that $(x,t)=1$ and therefore
$(y,t)=t^2$. Setting $x_1=y$ and $y_1=x$, one has
$t_1=(y_1,x_1)=t\inv$, $(y_1,t_1)=1$ and
$(x_1,t)=t_1^2$. Therefore $x_1, y_1$ and $t_1$
satisfy one of the conditions of 2 and this is
equivalent with $x, y$ and $t$ satisfying one of
the conditions of 2'.
\end{proof}

We will need the following remark.

\begin{remark}\label{t4}
{\rm It follows from the proof of
Lemma~\ref{TodosLosCasosRango2} that if $G$ is a
non-abelian $2$-group of Kleinian type with a
commutator $t$ of order $4$ then there exist
$x_1,y_1,x_2,y_2,x'_2,y'_2\in G$ with
$t=(y_1,x_1)=(y_2,x_2)=(y'_2,x'_2)$, and so that
$x_1$ and $y_1$ satisfy condition 1, $x_2$ and
$y_2$ satisfy condition 2, and $x'_2$ and $y'_2$
satisfy condition 2'. Moreover if $x_1$ and $y_1$
satisfy 1(a) (respectively, 1(b) or 1(c)) then
$x_2$ and $y_2$ satisfy 2(a) (respectively, 2(b) or
2(c)) and $x'_2$ and $y'_2$ satisfy 2'(a)
(respectively, 2'(b) or 2'(c)).}
\end{remark}

\begin{lemma}\label{lema6}
Let $G=\GEN{x_{1},x_{2},x_{3}}$ be a non-abelian
$2$-group such that  $x_{2}^{4}\neq 1$, $x_{3}\in
Z(G)$, $(x_{1},t)=t^{2}\neq 1$ and $(x_{2},t)=1$
with $t=(x_{2},x_{1})$. If $G$ is of Kleinian type
then $x_{3}^{2}\in Z(\GEN{x_{1},x_{2}})^{2}$.
\end{lemma}
\begin{proof}
We argue by contradiction. So suppose
$G=\GEN{x_1,x_2,x_3}$ has minimal order among the
possible counterexamples to the lemma. In particular
$x_3^2\ne 1$. By statement~\ref{ExpCentro} of
Lemma~\ref{SSP-index}, $x_3$ has order $4$. Let
$G_{1}=\GEN{x_{1},x_{2}}$. By
Lemma~\ref{TodosLosCasosRango2},
$Z(G_1)=\GEN{x_1^2,t^2,tx_2^2}$ and therefore
$Z(G_1)^2=\GEN{x_1^4,t^2x_2^4}$.

Suppose $z\in Z(G_1)$ is such that $1\neq z$,
$t^{2}\not\in \GEN{z}$ and $x_{2}^{4}\not\in
\GEN{z}$. Then $G_{1}'\cap \GEN{z}=\GEN{t}\cap
\GEN{z}=1$ and the minimality of the order of $|G|$
applied to the group $G/\GEN{z}$ yields
$x_{3}^{2}\in \GEN{z,Z(G_{1})^{2}}$.

If $t^{2}\neq x_{2}^{4}$ then $z=t^{3}x_{2}^{2}$ is
a non trivial central element of order $4$ so that
$t^{2}\not\in \GEN{z}$ and $x_{2}^{4}\not\in
\GEN{z}$. Hence, by the previous, $x_3^2$ is an
element of order $2$ of the group
$H=\GEN{t^{3}x_{2}^{2}, Z(G_{1})^{2}}$. Since
$t^3x_2^2\in Z(G_1)$, one has that $H=Z(G_1)^2 \cup
t^3x_2^2 Z(G_1)^2$. Then $x_3^2 = t^3 x_2^2 w^2$ for
some $w\in Z(G_1)$. So $1=x_3^4 = (t^3 x_2^2)^2 =
t^{2}x_{2}^{4} \ne 1$, a contradiction.

Thus we have that $t^{2}=x_{2}^{4}\ne 1$.
Lemma~\ref{TodosLosCasosRango2}  therefore implies
that we have one of the following properties: (i)
$x_{1}^{2}=t^{2}=x_{2}^{4}$ or $t=x_{1}^{\pm
2}x_{2}^{-2}$  (this is case 2(b)), or (ii)
$x_{1}^{4}=1$ and $t=x_{2}^{-2}$ (this is case
2(c)). In both cases we have $x_{1}^{4}=1$ and
therefore $Z(G_{1})^{2}=1$. Thus, if $z\in Z(G_{1})$
has order $2$ and $z\neq t^{2}=x_{2}^{4}$ then, by
the above argument, we have that $x_{3}^{2}\in
\GEN{z,Z(G_{1})^{2}}=\GEN{z}$; hence $x_{3}^{2}=z$.
This shows that $x_3^2$ is the unique central
element of order $2$ in $Z(G)$. Therefore $Z(G)$ is
cyclic generated by $x_3$, $Z(G_1)=\GEN{t^2}$ and
$x_3^3 = t^2$. Since $tx_{2}^{2}$ is central of
order  at most $2$ we thus get that either
$tx_{2}^{2}=1$ or $tx_{2}^{2}=t^{2}$, that is,
$t=x_{2}^{\pm 2}$.

Then $K=\GEN{x_{2},x_{3}}$ is an abelian subgroup of
index $2$ in $G$. Let $H=\GEN{tx_3^{-1}}$. Clearly
$K/H$ is cyclic (generated by $\overline{x_2}$).
Thus  $K=N_{G}(H)$  and $(K,H)$ is a strong Shoda
pair of $G$.
Using also statement~\ref{Indices} of
Lemma~\ref{SSP-index},  it follows that $[K:H]\leq
4$ and hence $t^{2}=x_{2}^{4}\in H=\{1,
tx_{3}^{-1}\}$. Then $t=x_3\inv \in Z(G)$, a
contradiction.
\end{proof}

\begin{lemma}\label{ComCic}
Let $G=\GEN{x_{1},x_{2},x_{3}}$ be a $2$-group of
Kleinian type with $G'=\GEN{t}$ of order $4$. Let
$t_{ij}=(x_{j},x_{i})$ with $1\leq i<j\leq 3$.
Assume that $t=t_{12}$, $(x_{1},t)=t^{2}$,
$(x_{2},t)=1$ and $t_{23}=1$.
\begin{enumerate}
\item If $t_{13}\in\GEN{t^{2}}$ then $x_{3}^{4}=1$. If, moreover, $x_{2}^{4}\neq 1$ then either $t_{13}=1$ and
$x_{3}^{2}\in Z(\GEN{x_{1},x_{2}})^{2}$ or
$t_{13}=t^2$ and $x_2^4x_3^2 \in
Z(\GEN{x_{1},x_{2}})^{2}$.
\item If $t_{13}\not\in\GEN{t^{2}}$ then $x_{3}^{4}=x_{2}^{4}$.
If, moreover, $x_{2}^{4}\neq 1$ then either
$t_{13}=t^{-1}$ and $x_{2}^{2}x_{3}^{2}\in
Z(\GEN{x_{1},x_{2}})^{2}$ or $t_{13}=t$ and
$t^{2}x_{2}^{2}x_{3}^{2}\in
Z(\GEN{x_{1},x_{2}})^{2}$.
\end{enumerate}
\end{lemma}

\begin{proof}
1. Assume first that $t_{13}\in\GEN{t^{2}}$. If
$t_{13}=1$ then $x_{3}\in Z(G)$ and if
$t_{13}=t^{2}$ then $x_{2}^{2}x_{3}\in Z(G)$. In
both cases, because of statement~\ref{ExpCentro} of
Lemma~\ref{SSP-index} and statement~\ref{Exponente}
of Lemma~\ref{PropEl}, we obtain that $x_{3}^{4}=1$.
The second statement is a consequence of
Lemma~\ref{lema6}, applied to $\GEN{x_1,x_2,x_3}$ if
$t_{13}=1$, and to $\GEN{x_1,x_2,x_2^2x_3}$ if
$t_{13}=t^2$.

2. Assume second that $t_{13}\not\in \GEN{t^{2}}$.
Then $t_{13}\in\{ t,t^{-1} \}$. If $t_{13}=t^{-1}$
then $x_{2}x_{3}\in Z(G)$ and hence
$(x_{2}x_{3})^{4}=x_{2}^{4}x_{3}^{4}=1$. So
$x_{2}^{4}=x_{3}^{4}$. If, moreover, $x_{2}^{4}\neq
1$ then by Lemma~\ref{lema6} we have that $x_{2}^{2}
x_{3}^{2}\in Z(\GEN{x_{1},x_{2}})^{2}$. On the other
hand if $t_{13}=t$ then $t x_{2} x_{3}\in Z(G)$ and
hence $(t x_{2}x_{3})^{4}=x_{2}^{4}x_{3}^4=1$. So
$x_{2}^{4}=x_{3}^{4}$. If, moreover, $x_{2}^{4}\neq
1$, then again by Lemma~\ref{lema6}, we have that
$t^{2} x_{2}^{2} x_{3}^{2}\in
Z(\GEN{x_{1},x_{2}})^{2}$. This finishes the proof.
\end{proof}

We need one more lemma before giving the proof of
(D) implies (F) for nilpotent groups.

\begin{lemma}\label{PrimerCaso}
Let $G$ be a finite non-abelian $2$-group of
Kleinian type. Assume $G'\subseteq Z(G)$ and
$Z(G/T)$ has exponent $2$ for every proper subgroup
$T$ of $G'$. Then $G$ is an epimorphic image of
either $C_2^n\times \W$, $\W_{1n}$ or $\W_{2n}$ for
some $n$.
\end{lemma}

\begin{proof}
Applying the assumptions for $T=1$ one deduces that
$Z(G)$ and $G'$ have exponent $2$. Then $G/Z(G)$ has
exponent $2$ and therefore $G$ has exponent $4$.

First, we prove that $G$ has an abelian subgroup of
index 2. Otherwise, by statement~\ref{TAmitsur} of
Lemma~\ref{SSP-index},
$G/Z(G)=\GEN{\overline{x_1}}_2 \times
\GEN{\overline{x_2}}_2 \times
\GEN{\overline{x_3}}_2$ for some $x_1,x_2,x_3\in G$.
Hence $G=\GEN{Z(G),x_{1},x_{2},x_{3}}$ and
$G'=\GEN{t_{12}}_2 \times \GEN{t_{13}}_2 \times
\GEN{t_{23}}_2$, where $t_{i,j}=(x_j,x_i)$ for
$1\leq i<j\leq 3$.  If $x\in G$ then
$G_x=\GEN{(x,y): y\in G}$ is a proper subgroup of
$G'$ and the image of $x$ in $G/G_x$ is central.
Therefore $x^2\in G_x$, by assumption. In
particular, $x_1^2=t_{12}^{\alpha_2}
t_{13}^{\alpha_3}$, $x_2^2=t_{12}^{\beta_1}
t_{23}^{\beta_3}$ and $x_3^2=t_{13}^{\gamma_1}
t_{23}^{\gamma_2}$, for some $\alpha_2$, $\alpha_3$,
$\beta_1$, $\beta_3$, $\gamma_1$, $\gamma_2 \in
\{0,1\}$. Then $(x_1x_2)^2=
t_{12}^{1+\alpha_2+\beta_1}
t_{13}^{\alpha_3}t_{23}^{\beta_3}
  \in\GEN{t_{12},t_{13}t_{23}}$,
$(x_1x_3)^2 = t_{12}^{\alpha_2}
t_{13}^{1+\alpha_3+\gamma_1}t_{23}^{\gamma_2}
  \in\GEN{t_{13},t_{12}t_{23}}$ and
$(x_2x_3)^2 = t_{12}^{\beta_1}
t_{13}^{\gamma_1}t_{23}^{1+\beta_3+\gamma_2}
  \in\GEN{t_{23},t_{12}t_{13}}$.
This implies that the $\alpha$'s, $\beta$'s and
$\gamma$'s with the same subindex are equal, so
$x_1^2=t_{12}^{a_2} t_{13}^{a_3}$,
$x_2^2=t_{12}^{a_1} t_{23}^{a_3}$ and
$x_3^2=t_{13}^{a_1} t_{23}^{a_2}$, for some
$a_1,a_2,a_3\in \{0,1\}$. Applying once more the
property to $x_1x_2x_3$ one obtains that
    $$t_{12}^{1+a_1+a_2} t_{13}^{1+a_1+a_3} t_{23}^{1+a_2+a_3} = t_{12} t_{13} t_{23} x_1^2 x_2^2 x_3^2=
    (x_1 x_2 x_3)^2 \in \GEN{t_{12}t_{13},t_{12}t_{23}}$$
and therefore $3+2a_1+2a_2+2a_3 \equiv 0 \mod 2$, a
contradiction.

Therefore $G=\GEN{x,y_1,\ldots,y_n}$ where
$Y=\GEN{y_1,\ldots,y_n}$ is an abelian subgroup of
index $2$ in $G$. In particular $G'\subseteq \langle
y_{1},\ldots , y_{n}\rangle$. If $y_{i}^2=1$ for
every $i=1,\dots,n$, then $G$ is an epimorphic image
of $\W_{1n}$, as desired. Otherwise, we may assume
without loss of generality that $y_1$ has exponent
$4$ and so $(y_1,x)\ne 1$, because $Z(G)$ has
exponent 2. If $|G'|=2$ then $(y_i,x)\in
\GEN{(y_1,x)}$ and, replacing $y_i$ by $y_1y_i$ if
needed, one may assume that $y_i\in Z(G)$ for every
$i\ge 2$. Then $G$ is a quotient of $\W\times
C_2^{n-1}$. Finally suppose that $|G'|>2$. Then, for
every $i\ge 2$, replacing $y_{i}$ by $y_{1}y_{i}$ if
needed, we may assume that $y_{i}$ has order 4. Then
$G_{y_i}=\GEN{(y_i,x)}$ is a proper subgroup of $G'$
and therefore $1\ne y_i^2 \in G_{y_i} =
\GEN{(y_i,x)}$ and so $y_i^2 = (y_i,x)$. It follows
that $G$ is an epimorphic image of $\W_{2n}$.
\end{proof}

We are ready to prove (D) implies (F) for nilpotent
groups. So let $G$ be a non-abelian finite nilpotent
group of Kleinian type. Hence by
statements~\ref{Descomposicion} and \ref{ExpCentro}
of Lemma~\ref{SSP-index} $G=G_3\times G_2$, where
$G_3$ is an elementary abelian 3-group, $G_2$ is a
non-abelian $2$-group and the exponent of
$Z(G)=G_3\times Z(G_2)$ divides $4$ or $6$.

We will deal separately with three cases. (1)
$G_3\ne 1$, (2) $G_3=1$ and $G'\subseteq Z(G)$ and
(3) $G_3=1$ and $G'\not\subseteq Z(G)$.

(1) Assume $G_{3}$ is not trivial. We will show that
$G_2$ satisfies the hypothesis of
Lemma~\ref{PrimerCaso}. This implies that $G_2$ is
isomorphic to a quotient of either $\W\times C_2^n$,
$\W_{1n}$ or $\W_{2n}$, for some $n$. Hence
condition (F.1) of Theorem~\ref{Main} holds.

If $T$ is a proper subgroup of $G_2'$ then, since
also $G/T$ is of Kleinian type, the exponent of
$Z(G/T)$ is $6$, by statement~\ref{ExpCentro} of
Lemma~\ref{SSP-index}. Hence $Z(G_2/T)$ has exponent
$2$, as desired.

Next we need to show that $G'\subseteq Z(G)$. We
prove this by contradiction. So assume that
$G'\not\subseteq Z(G)$. Then, by
statement~\ref{Derivado} of Lemma~\ref{SSP-index}
and statement~\ref{Derivado1} of Lemma~\ref{PropEl},
there exist $x,y\in G_2$ such that $t=(y,x)$ has
order $4$. Because of Remark~\ref{t4} one may assume
without loss of generality that $x$ and $y$ satisfy
condition 1 of Lemma~\ref{TodosLosCasosRango2}. So
$x^2,y^2 \in Z(G_2)$ and therefore $x^4=y^4=1$.
Since $t^2\ne 1$, case 1(a) does not hold and so
$t^2\in \{x^2,y^2,x^2y^2\}$. By symmetry one may
assume that $t^2=x^2$ or $t^2=x^2y^2$. Notice that
$(xy)^2 = t\inv x^2 y^2$ and therefore $xy$ has
order $8$. Thus $H=\GEN{x,y}$ is a non-abelian group
of exponent 8 which is an epimorphic image of one of
the following two groups:
    $$\matriz{{rcl} H_1 &=& \GEN{a,b| a^4 = b^4 = 1, t=(b,a) , t^a = t^b = t\inv, t^2=a^2} \\
    H_2 &=& \GEN{a,b| a^4 = b^4 = 1, t=(b,a) , t^a = t^b = t\inv, t^2=a^2b^2}}$$
On the other hand, $\Q(C_3\times Q_{16})$ has an
epimorphic image isomorphic to
$\Q(\xi_3)\otimes_{\Q} \HQ(\Q(\sqrt{2})) \cong
M_2(\Q(\xi_3,\sqrt{2}))$ and $\Q(C_3\times
D_{16}^-)$ has an epimorphic image isomorphic to
$\Q(\xi_3)\otimes_{\Q} M_2(\Q(\sqrt{-2})) \cong
M_2(\Q(\xi_3,\sqrt{-2}))$ (see (\ref{WDDQ})). Then,
statement~\ref{Indices} of Lemma~\ref{SSP-index}
implies that neither $C_3\times Q_{16}$ nor $C_3
\times D_{16}^-$ are of Kleinian type. Since
$H_1/\GEN{a^2b^2}\cong Q_{16}$ and
$H_2/\GEN{b^2}\cong D_{16}^-$, Lemma~\ref{Closed}
implies that neither $C_3\times H_1$ nor $C_3 \times
H_2$ are of Kleinian type. Since $|H|\ge 16$ and
$|H_1|=|H_2|=32$, we have that $H$ is a non-abelian
group of order $16$ with an element of order $8$.
This implies that $H$ is isomorphic to either
$D_{16}, D_{16}^+, D_{16}^-$ or $Q_{16}$. However
$H$ is not isomorphic to $D_{16}$ because the latter
is not of Kleinian type, and it is also not
isomorphic to $D_{16}^+$ because the commutator of
$D_{16}^+$ has order $2$. Moreover the same argument
as above shows that $H$ is not isomorphic to neither
$Q_{16}$ nor $D_{16}^-$. This yields in all cases a
contradiction. So $G'\subseteq Z(G)$ and this
finishes the proof of (1).

(2) Assume that $G_3=1$ and $G'\subseteq Z(G)$. We
prove that $G$ is isomorphic to a quotient of either
$\V\times A$, $\U_1\times A$, $\U_2\times A$,
$\V_{1n}$ or $\V_{2n}$, for  an abelian group $A$ of
exponent $4$. Hence condition (F.2) of
Theorem~\ref{Main} holds.

From statement~\ref{Exponente} of Lemma~\ref{PropEl}
and statement~\ref{ExpCentro} of
Lemma~\ref{SSP-index},  we know that the exponent of
$G$ divides $8$ and the exponent of $Z(G)$ divides
$4$. Moreover, the assumptions and
statement~\ref{Derivado} of Lemma~\ref{SSP-index}
imply that $G'$ is of exponent $2$. Hence $g^{2}\in
Z(G)$ for all $g\in G$.

If $G$ does not contain an abelian subgroup of index
$2$ then, by statement~\ref{TAmitsur} of
Lemma~\ref{SSP-index}, $G=\langle
Z(G),y_{1},y_{2},y_{3}\rangle$, and $G/Z(G)$ and
$G'=\GEN{t_{ij}=(y_{j},y_{i})\mid1\leq i<j\leq 3}$
are both elementary abelian groups of order $8$. By
statement~\ref{Derivado1} of Lemma~\ref{PropEl}, it
follows that there exist $\alpha_2$, $\alpha_3$,
$\beta_1$, $\beta_3$, $\gamma_1$ and $\gamma_2$ in
$\{0,1\}$ such that
$y_1^4=t_{12}^{\alpha_2}t_{13}^{\alpha_3}$,
$y_2^4=t_{12}^{\beta_1}t_{23}^{\beta_3}$  and
$y_3^4=t_{13}^{\gamma_1}t_{23}^{\gamma_2}$. Applying
again part~\ref{Derivado1} of Lemma~\ref{PropEl}, it
follows that $(y_1y_2)^4=t_{12}^{\alpha_2+\beta_1}
t_{13}^{\alpha_3}t_{23}^{\beta_3}
  \in\GEN{t_{12},t_{13}t_{23}}$,
$(y_1y_3)^4=t_{12}^{\alpha_2}
t_{13}^{\alpha_3+\gamma_1}t_{23}^{\gamma_2}
  \in\GEN{t_{13},t_{12}t_{23}}$ and
$(y_2y_3)^4=t_{12}^{\beta_1}
t_{13}^{\gamma_1}t_{23}^{\beta_3+\gamma_2}
  \in\GEN{t_{23},t_{12}t_{13}}$.
Hence $\alpha_3=\beta_3$, $\alpha_2=\gamma_2$ and
$\beta_1=\gamma_1$. To simplify notation,  put
$a_1=\beta_1$, $a_2=\alpha_2$ and $a_3=\alpha_3$.
Then, once more applying statement~\ref{Derivado1}
of Lemma~\ref{PropEl}, we get
    \begin{equation}\label{ocho}
    \matriz{{cc}
  \matriz{{rcl}
  y_1^4&=&t_{12}^{a_2}t_{13}^{a_3} \\
  y_2^4&=&t_{12}^{a_1}t_{23}^{a_3} \\
  y_3^4&=&t_{13}^{a_1}t_{23}^{a_2}} \quad \mbox{ and } \quad
  \matriz{{rcl}
  (y_1y_2)^4&=&t_{12}^{a_1+a_2}t_{13}^{a_3}t_{23}^{a_3} \\
  (y_1y_3)^4&=&t_{12}^{a_2}t_{13}^{a_1+a_3}t_{23}^{a_2} \\
  (y_2y_3)^4&=&t_{12}^{a_1}t_{13}^{a_1}t_{23}^{a_2+a_3}\\
  (y_1y_2y_3)^4&=&t_{12}^{a_1+a_2}t_{13}^{a_1+a_3}t_{23}^{a_2+a_3}.}}
  \end{equation}
Because each $a_{i}\in \{ 0,1 \}$, we obtain that at
least one of the seven elements in (\ref{ocho}) is
equal to $1$. Without loss of generality, we may
assume that $y_1^4=1$, and hence $a_2=a_3=0$. Then
$y_2^4=t_{12}^{a_1}$ and $y_3^4=t_{13}^{a_1}$. If
$a_{1}=0$ then it follows that $G$ is an epimorphic
image of $\U_{1}\times C_{4}^{n}$ for some $n$. If
$a_{1}=1$ then $G$ is an epimorphic image of
$\U_{2}\times C_{4}^{n}$ for some $n$.

We now consider the case that $G$ has an abelian
subgroup $\langle y_{1},\ldots , y_{n} \rangle$ of
index $2$. Write $G=\langle x,y_{1},\ldots , y_{n}
\rangle$. If $y_i^4=1$ for every $i$, then $G$ is an
epimorphic image of $\V_{1n}$. So assume that some
$y_i$, say $y_1$, has order 8. In particular
$(y_1,x)\ne 1$. As in the case where $G_3\ne 1$ and
$G_2$ has an abelian subgroup of index 2, if
$|G'|=2$ then one may assume that $y_i$ is central
for every $i\ge 2$ and therefore $y_i^4=1$. This
implies that $G$ is an epimorphic image of $\V\times
C_4^{n-1}$. Finally, assume that $y_1$ has order 8
and $|G'|>2$. Again following the same pattern as in
the case of $G_3\ne 1$, replacing $y_{i}$ by
$y_{1}y_{i}$ one may assume that each $y_{i}$ has
order $8$ and applying statement~\ref{Derivado1} of
Lemma~\ref{PropEl}, one deduces that $y_i^4=(y_i,x)$
for every $i$. It follows that $G$ is an epimorphic
image of $\V_{2n}$.

(3) Assume that $G_3=1$ and $G'\not\subseteq Z(G)$.
We prove that $G$ is an epimorphic image of either
$\T\times A$, $\T_{1n}$, $\T_{2n}$ or $\T_{3n}$ for
$A$ an elementary abelian $2$-group.

By statement~\ref{Exponente} of Lemma~\ref{PropEl},
the exponent of $G$ divides $8$. By
statements~\ref{TAmitsur} and \ref{Derivado} of
Lemma~\ref{SSP-index}, $G$ has an abelian subgroup
$Y=\langle y_{1},\ldots , y_{n} \rangle$ of index
$2$ and $G'$ has exponent $4$. Then $G=\GEN{Y,x}$
for some $x\in G$ and $G'=\GEN{t_1,\dots,t_n}$,
where $t_i=(y_i,x)$ for $i=1,\dots,n$. We may
assume, without loss of generality, that $t_{1}$ is
of order $4$ (and thus $t_{1}$ is not central).
Since $G'\subseteq Y$, $(t_{i},y_{j})=1$ for all
$1\leq i,j\leq n$. If $t_{j}$ is not central then,
by Lemma~\ref{TodosLosCasosRango2},
$(x,t_{j})=t_{j}^{2}$. If $t_{j}$ is central, and
thus of order two, we also get that
$(x,t_{j})=1=t_{j}^{2}$. So, in all cases we have
$(x,t_{j})=t_{j}^{2}$.

We now show that we may assume that $\GEN{t_{1}}\cap
\GEN{t_{i}}=1$ for every $i\geq 2$.

Because the order of $t_{i}$ divides $4$, this is
clear if $t_{i}^{2}\not\in \GEN{t_{1}}$. If
$t_{i}\in \GEN{t_{1}}$, say $t_i=t_{1}^{a}$ then we
replace $y_{i}$ by $y'_{i}=y_{1}^{-a} y_{i}$ to make
$(y_{i}',x)=1$, because $(y_{i}',x)=(y_{1}^{-a}
y_{i},x) = (y_{i},x)^{y_{1}^{a}} (y_{1}^{-a},x) =
t_{i} t_{1}^{-a} = 1$. In the remaining case
$t_{i}\not\in \GEN{t_{1}}$ and
$t_{i}^{2}\in\GEN{t_{1}^{2}}$. Then either
$t_{i}^{2}=1$ or $t_{i}^{2} = t_{1}^{2}$. If
$t_{i}^{2} = 1$ then the claim is clear. If
$t_{i}^{2}=t_{1}^{2}$ then replacing $y_{i}$ by
$y_{i}'=y_{1}y_{i}$ we obtain that $ (y_{i}',x) =
(y_{1}y_{i},x) = t_{i}^{y_{1}\inv} t_{1} =t_{i}t_{1}
\not\in \GEN{t_{1}} $ and $(y_{i}',x)^{2}=1$, which
finishes the proof of the claim. So from now on we
assume that for $i\geq 2$, $\langle t_{1} \rangle
\cap \langle t_{i} \rangle =1$. Since the order of
$t_i$ divides $4$ and the order of $t_1$ is $4$,
this implies that $\GEN{t_i}\cap \GEN{t_1t_i}=1$ for
$i\ge 2$.

For $i=1,\ldots,n$ we put
$F_{i}=\GEN{x,y_{1},y_{i}}$ and we prove three
claims.

{\it Claim 1}: If $y_1^{4}=1$ then $y_{j}^{4} = 1$
for every $j$ (with $1\leq j \leq n$).

Indeed, suppose $y_{1}^{4}=1$. If $t_{j}=1$ then
$y_{j}$ is central in $G$ and thus, by
statement~\ref{ExpCentro} of Lemma~\ref{SSP-index},
we get at once that $y_{j}^{4}=1$. So assume that
$t_j\ne 1$. We now apply statement~\ref{Derivado1}
of Lemma~\ref{PropEl} to the group $F_j$. Since
$\GEN{t_1}\cap \GEN{t_j}=1$, one has $(F_j)_{y_j} =
\GEN{t_j} \ne F'_j \ne \GEN{y_1y_j} = (F_j)_{y_1y_j}
$ and hence $y_j^4\in \GEN{t_j}$ and $y_j^4 =
(t_1t_j)^4 \in \GEN{t_1t_j}$. Thus $y_{j}^{4} \in
\GEN{t_{j}}\cap \GEN{t_{1}t_{j}}=1$. This proves
Claim 1.

{\it Claim 2}: If $x^{4}=1$ and $t_1=y_1^{-2}$ then
$t_{j}=y_{j}^{-2}$ for every $j$ (with $1\leq j \leq
n$).

Indeed, suppose $x^4=1$ and $t_{1}=y_{1}^{-2}$. Let
$Z=\GEN{x^2,t_1^2}$. Then $Z$ is a subgroup of
$Z(F_1)$. Moreover $y_1^{-2} = t_1\not\in Z(F_1)$
and $(xy_1^i,t_1)=t_1^2\ne 1$ for every $i$. This
shows that $Z=Z(F_1)$. Hence $Z(F_1)^2=1$.

Let $j$ be such that $2\leq j \leq n$. Since
$t_{1}^{2}\not\in \GEN{t_{j}}$ (because
$\GEN{t_{1}}\cap \GEN{t_{j}} =1$ and $t_{1}$ has
order $4$) we can apply Lemma~\ref{lema6} to the
elements $x_1=\overline{x}$, $x_2=\overline{y_1}$
and $x_3=\overline{y_j}$ of the non-abelian Kleinian
group $F_{j}/\GEN{t_{j}}$ and deduce that
$\overline{y_{j}}^{2}\in
\overline{Z}^{2}=1$. 
Hence
\begin{equation} \label{eqclaim2}
y_{j}^{2}\in \GEN{t_{j}}
\end{equation}

We now proceed by considering the possible orders of
$t_{j}$. If $t_{j}=1$ then (\ref{eqclaim2}) implies
that $y_{j}^{-2}=1=t_{j}$, as desired. If $t_{j}$
has order $4$ then, again because $\GEN{t_{1}}\cap
\GEN{t_{j}}=1$, the second part of
Lemma~\ref{ComCic} is applicable to the group
$F_{j}/\GEN{t_{1}^{-1} t_{j}}$,  for
$x_1=\overline{x}$, $x_2=\overline{y_1}$ and
$x_3=\overline{y_j}$.  It follows that
$\overline{t_{1}^{2}y_{1}^{2}y_{j}^{2}}\in
\overline{Z}^{2}=1$. Hence
$t_{1}^{2}y_{1}^{2}y_{j}^{2}=t_{1}y_{j}^{2}\in
\GEN{t_{1}^{-1}t_{j}}$. Combining this with
(\ref{eqclaim2}), we obtain $y_{j}^{2}\in
\GEN{t_{j}}\cap t_{1}^{-1}\GEN{t_{1}^{-1}t_{j}} =
\{t_{j}^{-1}\}$, as desired. If $t_{j}$ has order
two then again (\ref{eqclaim2}) implies that either
$y_{j}^{2}=t_{j}$ or $y_{j}^{2}=1$. The former is as
desired. In the second case we can apply
Lemma~\ref{lema6} to the non-abelian Kleinian group
$F_{j}/\GEN{t_{1}^{2}t_{j}}$  (note that
$\overline{t_{1}y_{j}}$ is central in this group).
It follows that $(\overline{t_{1}y_{j}})^{2} \in
\overline{Z^{2}}=1$ and thus  $t_{1}^{2} =
t_{1}^{2}y_{j}^{2} \in \GEN{t_{1}^{2}t_{j}}$. Hence
$t_{1}^{2}=t_{1}^{2}t_{j}$, a contradiction. This
finishes the proof of Claim 2.

{\it Claim 3}: If $G'$ is not cyclic then
$y_{i}^{4}\in \GEN{t_{i}}$ for every $i$ with $1\leq
i \leq n$. If, furthermore, $t_{i}^{2}=y_{i}^{4}\neq
1$ for some $i\geq 1$, then $x^{4}=1$.

Assume that $G'$ is not cyclic. Then
$G_{y_i}=\GEN{t_i}\ne G'$, for each $i\ge 1$. Hence,
by statement~\ref{Derivado1} of Lemma~\ref{PropEl},
$y_{i}^{4} \in \GEN{t_{i}}$, as desired. Assume,
furthermore, that $x^{4}\neq 1$ and
$t_{i}^{2}=y_{i}^{4}\neq 1$ for some $i\geq 1$. By
Lemma~\ref{TodosLosCasosRango2},
$t_{i}y_{i}^{2}=x^{\pm 2}$ and therefore
$1=t_{i}^{2}y_{i}^{4}=x^{4}$, a contradiction. Hence
the claim follows.

We now consider 3 cases.

{\it Case 1}. Suppose $y_{1}^{4}=1$.

Because of Claim 1 we obtain that $y_{j}^{4}=1$ for
every $j$. Hence we conclude that $G$ is a quotient
of $\T_{1n}$.

{\it Case 2}. Suppose $x^{4}=1$ and
$t_{1}=y_{1}^{-2}$.

Because of  Claim 2 we conclude that $G$ is a
quotient of $\T_{2n}$.

{\it Case 3}. Suppose that neither Case 1 nor Case 2
hold.

{\it Claim 4}. One may assume that, for every $i\ge
1$, if $t_{i}^{2}\neq 1$ then one has $y_{i}^{4}
\neq 1$, $t_{i} \neq y_{i}^{-2}$ and either
$t_iy_i^2=x^{\pm 2}$ or $x^2 \in \{t_i^2,y_i^4\}$.

Suppose that $t_i^2\ne 1$. Then $x$ and $y_i$
satisfy condition 2 of
Lemma~\ref{TodosLosCasosRango2} and hence one of the
three cases (a), (b) or (c) of this statement holds.
If $y_i^4 = 1$ then interchanging the roles of $y_1$
and $y_i$ one may assume that Case 1 holds. So one
may assume that $y_i^4\ne 1$ and hence (a) does not
hold. Suppose now that $t_i=y_i^{-2}$. If $x^4 \ne
1$ then (c) does not hold and from (b) we get that
$1=t_{i}y_{i}^{2}=x^{\pm 2}\neq 1$, a contradiction.
Thus in this case $x^4 = 1$ and interchanging again
the roles of $y_1$ and $y_i$ one may assume that
Case 2 holds. So one may assume that $y_i^4\ne 1$
and $t_i \ne y_i^{-2}$. Then neither (a) nor (c)
holds. Thus (b) holds and this finishes the proof of
Claim 4.

Assume that $G'$ is cyclic. Thus
$\GEN{t_{i}}\subseteq \GEN{t_{1}}$. Since we already
know that $\GEN{t_{1}}\cap \GEN{t_{i}} =1$, for
$i\geq 2$, this implies that $t_{i}=1$. Hence
$y_{i}$ is central for $i\ge 2$. By
Lemma~\ref{lema6}, we get that $y_{i}^{2}=z_{i}^{2}$
for some $z_i\in Z(F_{1})$. Then, replacing $y_{i}$
by $y_{i}z_{i}$, we may assume that $y_{i}^{2}=1$,
for $i\ge 2$. Thus $G$ is an epimorphic image of
$F_1\times C_2^{n-1}$. We are going to show that
$F_1$ is an epimorphic image of $\T$ and therefore
$G$ is an epimorphic image of $\T\times C_2^{n-1}$.

First assume that $x^4\ne 1$ and so, by Claim 4,
$t_1y_i^2 = x^{\pm 2}$. If $t_1y_i^2 = x^{-2}$ then
$(xy_1)^2 = t_1\inv x^2 y_1^2 = x^4 y_1^4 = t_1^2 =
(xy_1,t_1)$ and $t_1=(y_1,xy_1)$. Then replacing $x$
by $xy_1$ one sees that $F_1$ is an epimorphic image
of $\T$. A similar computation shows that if
$t_1y_1^2 = x^2$ then replacing $x$ by $xy_1\inv$
one deduces that $F_1$ is an epimorphic image of
$\T$. Second assume that $x^{4}=1$. By Claim 4,
either $t_{1}y_{1}^{2}=x^{2}$, $x^{2}=y_{1}^{4}$ or
$x^{2}=t_{1}^{2}$. If $x^{2}=t_{1}^{2}$ then $H$ is
clearly an epimorphic image of $\T$. If
$x^{2}=y_{1}^{4}$ then $x\mapsto xy_1^2$ and
$y\mapsto y_1$ induces an epimorphism $\T\rightarrow
H$ and if $t_{1}y_{1}^{2}=x^{2}$ one gets an
epimorphism $\T \rightarrow H$ given by $x\mapsto
xy_1$ and $y\mapsto y_1$. This finishes the proof if
$G'$ is cyclic.

Assume that $G'$ is not cyclic.

By Claim 3, $y_{i}^{4}\in \GEN{t_{i}} \cap Z(G)$ for
every $i$. In particular, if $t_{i}^{2}\neq 1$ (for
example for $i=1$), then $y_{i}^{4}=t_{i}^{2}$
(because by assumption $y_{i}^{4}\neq 1$ by Claim
4). Because of Claim 3 we then get that $x^{4}=1$.
Moreover $Z(\GEN{x,y_i}) = \GEN{x^2,t_i^2,t_iy_i^2}$
(see Lemma~\ref{TodosLosCasosRango2}) and so
$Z(\GEN{x,y_i})^2=1$ for every $i$ such that
$t_i^2\ne 1$.

We claim that $y_{i}^{2}\in \GEN{t_i}$ and
$t_{i}^{2}=1$ for every $i\geq 2$. Clearly the image
$\overline{y_{i}}$ of $y_{i}$ in
$H=F_{i}/\GEN{t_{i}}$ is central. Because
$\GEN{t_{i}} \cap \GEN{t_{1}}= 1 $ and
$y_{1}^{4}=t_{1}^{2}$, Lemma~\ref{lema6} is
applicable to the group $H$. Indeed,
$H=\GEN{x_1=\overline{x},x_2=\overline{y_1},x_3=\overline{y_i}}$,
$x_3\in Z(H)$, $(x_1,t)=x_2^4=t^2\ne 1$ and
$(x_2,t)=1$, where $t=(x_2,x_1)=\overline{t_1}$,
because $t_1^2\not\in\GEN{t_i}$. Therefore we get
that $x_3^{2}\in Z(\GEN{x_1,x_2})^{2}=1$, or
equivalently $y_{i}^{2}\in \GEN{t_{i}}$, as desired.
Assume now that $t_{i}^{2}\neq 1$. Then, by the
previous paragraph, $y_i^4=t_i^2$ and hence
$y_{i}^{2} =t_{i}$, (because the option $y_{i}^{2}=
t_{i}^{-1}$ is excluded by Claim 4). The last part
of Claim 4 now implies $x^{2}=t_{i}^{2}$.
Interchanging the role of $y_{1}$ and $y_{i}$ in the
above reasoning, we get that $t_{1}^{2}=x^{2}$.
Hence $t_{1}^{2}=t_{i}^{2}\ne 1$, contradicting with
$\GEN{t_{1}}\cap \GEN{t_{i}} =1$. This proves that
$t_{i}^{2}=1$ and shows the claim.

Let $i\geq 2$. The natural image of $t_{1}y_{i}$ is
central in the non-abelian Kleinian group
$F_{i}/\GEN{t_{1}^{2}t_{i}}$. Hence applying
Lemma~\ref{lema6} to this group, we obtain that
$\overline{t_{1}^{2}} \overline{y_{i}^{2}}\in
Z(\GEN{ \overline{x},\overline{y_{1}} })^{2}=1$.
Consequently $t_{1}^{2}y_{i}^{2}\in
\GEN{t_{1}^{2}t_{i}}$. Thus $y_{i}^{2} \in \{ 1,
t_{i}\} \cap \{ t_{1}^{2} , t_{i}\} =\{ t_{i} \}$,
i.e. $y_{i}^{2}=t_i$. Moreover, Claim 4 implies that
either $x^{2}=t_{1}^{2}$ or $t_{1}y_{1}^{2}=x^{2}$.
In the first case, $G$ is a quotient of $\T_{3n}$.
In the second case, setting $x'=y_1x$ and
$y_1'=y_1$, one has $t_{1}' = (y_{1},y_{1}x) =
t_{1}$ and $t_{i}'=(y_{i},y_{1}x)=t_{i}$. So
$y_{i}^{2}=t_{i}'$ for every $i\geq 2$ and ${x'}^{2}
= y_{1}xy_{1}x = t_{1}xy_{1}t_{1}xy_1=
t_{1}xt_{1}y_{1}xy_{1} = xy_{1}xy_{1} =
xt_{1}xy_{1}^{2} = t_{1}^{3}x^{2} y_{1}^{2} =
t_{1}^{3}t_{1}y_{1}^{2}y_{1}^{2} =y_{1}^{4} =
t_{1}^{2} =(t_{1}')^{2}$. This implies that again
$G$ is a quotient of $\T_{3n}$. It also finishes the
proof of (D) implies (F) for nilpotent groups.
$\qed$

\section{(D) implies (F), for non-nilpotent groups}
\label{non-nilpotent-case}

In this section we prove that (D) implies (F) for
finite groups that are not nilpotent.

Let $G$ be a finite non-nilpotent group of Kleinian
type. By statement~\ref{Descomposicion} of
Lemma~\ref{SSP-index}, $G$ is a semidirect product
$G_3\rtimes G_2$ of an elementary abelian $3$-group
$G_3$ and a $2$-group $G_2$. Moreover, since by
assumption $G$ is not nilpotent,
statement~\ref{TAmitsur} of Lemma~\ref{SSP-index}
implies that $G$ has an abelian subgroup $G_3\times
N_2$ such that $N_2$ has index $2$ in $G_2$. Thus
$G_{2}=N_{2}\cup N_{2}x$, for every $x\in
G_2\setminus N_2$. Let $K=G_{3}\cap Z(G)$. Then
$G_{3}=K\times M$, for some subgroup $M$ of $G_{3}$.
Note that $M$ is not trivial because, by assumption,
$G$ is not nilpotent. For every $m\in M$, let $k_m=m
m^x$ and $\widetilde{m}=k_m m$. Since $x^2$
centralizes $G_3$ and $G_3$ is abelian, $k_m\in K$
and thus $\widetilde{m}^x = k_m m^x = k_m^2 m\inv =
\widetilde{m}\inv$. Furthermore $k_{m_1m_2} = m_1
m_2 (m_1 m_2)^x = k_{m_1} k_{m_2}$ and hence
$\widetilde{m_1 m_2} = k_{m_1m_2} m_1 m_2 = k_{m_1}
m_1 k_{m_2} m_2 = \widetilde{m_1} \widetilde{m_2}$.
Hence $\widetilde{M}=\{\widetilde{m}\mid m\in M\}$
also is an elementary abelian $3$-group and $G_3 = K
\times \widetilde{M}$. So, replacing $M$ by
$\widetilde{M}$ we may assume that
$\widetilde{M}=M$, $a^x = a$ if $a\in K$, and $a^x =
a\inv$ if $a\in M$. Consequently,
    $G=K\times (M \rtimes G_2)=
    K \times (M \rtimes\GEN{N_2,x}),$
where $K$ and $M$ are elementary abelian $3$-groups,
$G_2=\GEN{N_2,x}=N_{2}\cup N_{2}x$ is a $2$-group,
$\GEN{N_{2},M}=N_{2}\times M$ is abelian and $x$
acts on $M$ by inversion. Notice that this group is
completely determined by $N_2$, $G_2$ and the ranks
$k$ and $m$ of $K$ and $M$ respectively. To
emphasize this information we use the following
notation
\begin{equation}\label{NNilp}
\matriz{{c}G= G_{k,m,N_2, G_2} =K\times (M \rtimes
G_2)=
    K \times ((M \times N_2): \GEN{\overline{u}}_2 ), \\ \hspace{0.2cm}
    (k\geq 0, m\geq 1, K=C_3^k, M=C_3^m, u^2\in N_2 \mbox{ and } w^u =w\inv \mbox{ for } w\in M)}
\end{equation}

%

Since Theorem~\ref{Main} has already been proved for
nilpotent groups, if $G$ is of Kleinian type and
$G_2$ is non-abelian then $K\times G_2$ satisfies
either condition (F.1), (F.2) or (F.3). In
particular, if $K\ne 1$ and $G_2$ is non-abelian,
then $K\times G_2$ satisfies condition (F.1) and so
the exponent of $G_2$ is $4$, $G_2'\subseteq Z(G_2)$
and the exponent of $Z(G_2)$ is $2$. In the
following four lemmas we find more restrictions on
$k$, $m$, $N_2$ and $G_2$.

\begin{lemma}\label{C_8C_8}
If $G=G_{k,m,N_2,G_2}$ is of Kleinian type then the
exponent of $G_2$ divides $8$. Furthermore,   if
$k\neq 0$ then the exponent of $G_2$ divides  $4$
and the exponent of $G_2\cap Z(G)$ is $1$ or $2$.
\end{lemma}

\begin{proof}
First we prove that the exponent of $G_2$ divides
$8$. This is a consequence of part~\ref{Exponente}
of Lemma~\ref{PropEl}, if $G_2$ is non-abelian. If
$G_2$ is abelian and $g\in G_2$ then $g^2\in N_2$
and therefore it is central. By
statement~\ref{ExpCentro} of Lemma~\ref{SSP-index}
the order of $g^2$ divides $4$ and thus the order of
$g$ divides $8$.

Assume now that $k\ne 0$, or equivalently $K\ne 1$.
If $G_2$ is non-abelian then $G_2$ has exponent $4$
and $Z(G_2)$ has exponent $2$ and so $G_2\cap Z(G)$
has exponent $1$ or $2$ as wanted. Otherwise, that
is if $G_2$ is abelian, $N_2=G_2\cap Z(G)$. Since
$K$ has a central element of order $3$, the exponent
of $Z(G)$ is either $3$ or $6$ and therefore the
exponent of $G_2\cap Z(G)$ is either $1$ or $2$.
Furthermore $g^2 \in G_2\cap Z(G)$ for every $g\in
G_2$ and so $g^4=1$.
\end{proof}

 Notice that if $M_1$ is a maximal subgroup of $M$
then $G/(K\times M_1) \cong G_{0,1,N_2,G_2}$. So to
obtain restrictions on $G_2$ and $N_2$ one may
assume without loss of generality that $k=0$ and
$m=1$. This will be used in the proof of the next
three lemmas.

\begin{lemma}\label{D8}
Assume that $G=G_{k,m,N_2,G_2}$ is of Kleinian type.
If $L$ is a normal subgroup of $G_{2}$ contained in
$N_{2}$ such that $G_{2}/L\cong D_{8}$ and  $a\in
N_{2}$ then $a^{2}\in L$. In particular, if
$G_{2}=D_{8}$ and $a$ is an element of order $4$ in
$G_{2}$ then $a\not\in N_{2}$.
\end{lemma}

\begin{proof}
 One may suppose that $k=0$ and $m=1$. Then $G/L \cong C_3 \rtimes D_8$. If $a^{2}\not\in L$ and $a\in N_{2}$ then
$G/L = \GEN{c}_3 \rtimes (\GEN{\overline{a}}_4
\rtimes \GEN{b}_2) = \GEN{c\overline{a}}_{12}
\rtimes \GEN{b}_2 \cong D_{24}$, contradicting
statement~\ref{Dihedral} of Lemma~\ref{SSP-index}.
\end{proof}

\begin{lemma}\label{Elementales}
Assume that $G=G_{k,m,N_2,G_2}$ is of Kleinian type.
Let $L$ be a normal subgroup of $G_{2}$ contained in
$N_{2}$. Then $G_{2}/L$ is not isomorphic to any of
the groups $Q_{16}$, $D_{16}^-$, $D_{16}^+$,
$\mathcal{D}$, $\mathcal{D}^{+}$.

In particular, if, moreover, $G_{2}/L$ is
non-abelian and has order $16$ then $G_2/L$ has
exponent $4$ and the exponent of $Z(G_2/L)$ is $2$.
\end{lemma}

\begin{proof}
We may assume that $k=0$ and $m=1$ and hence
$H=G/L=G_{0,1,Q,P}=\GEN{w}_3 \rtimes P$, where
$P=G_{2}/L$ and $Q=N_2/L$.

First assume that $P=Q_{16}=\GEN{a,b \mid
a^8=b^2a^4=1, ba=a^{-1}b}$ or $P=D_{16}^-=\GEN{a,b
\mid a^8=b^2=1,ba=a^3b}$. Then $P/\GEN{a^4}\cong
D_8$ and $a^2\not\in \GEN{a^4}$. By Lemma~\ref{D8},
$a\not\in Q$. However $(\GEN{w,a^2},1)$ is a strong
Shoda pair of $H$ and $[H: \GEN{w,a^2}]=4$,
contradicting statement~\ref{Indices}(a) of
Lemma~\ref{SSP-index}.

Second assume that $P=D_{16}^+ =\langle a,b \mid
a^{8}= b^{2}=1,\; ba=a^{5}b \rangle$. If $b\in Q$
then $a\not\in Q$ and $(A=\GEN{w,a^2,b},B=\GEN{b})$
is a strong Shoda pair of $G$ with $[A:B]=12$ and
$B$ is not normal in $G$, contradicting
statement~\ref{Indices}(b) of Lemma~\ref{SSP-index}.
On the other hand, if $b\not\in Q$ then,
interchanging generators if needed, we may assume
that $a\in Q$ and hence $(A=\GEN{w,a},1)$ is a
strong Shoda pair of $H$. Let $e=e(H,A,1)$, a
primitive central idempotent of $\Q H$. Then
$b^{2}e=e$ but $be\neq \pm e$ (because $be$ cannot
be central in $\Q He$). Hence $\Q H e$ is split.
Since $|A|=24$ we obtain a contradiction with
statement~\ref{Indices}(c) of Lemma~\ref{SSP-index}.

Third, assume  that $P=\mathcal{D}=\langle a,b,c
\mid ca=ac,\; cb=bc,\; a^{2}= b^{2}= c^{4}=1,\;
ba=c^{2}ab \rangle$. Since $ab$ is of order $4$ and
$\GEN{a,b}=D_{8}$,  Lemma~\ref{D8} implies that
$ab\not\in Q$. It thus follows that either $a\not\in
Q$ or $b\not\in Q$. By symmetry, we may assume that
$a\not\in Q$ and $b\in Q$. If $c\in Q$ then
$(A=\GEN{w,c,b},B=\GEN{b})$ is a strong Shoda pair
of $H$, $[A:B]=12$ and $\GEN{b}$ is not normal in
$H$, contradicting statement~\ref{Indices}(b) of
Lemma~\ref{SSP-index}. If $c\not\in Q$ then
$(A=\GEN{w,c^2,b},B=\GEN{b})$ is a strong Shoda pair
of $H$ such that $[H:A]=4$, again in contradiction
with statement~\ref{Indices}(b) of
Lemma~\ref{SSP-index}.

Fourth, assume that $P=\mathcal{D}^{+}=\langle a,b,c
\mid ca=ac,\; cb=bc,\; a^{4}=b^{2}=c^{4}=1,\;
ba=ca^{3}b \rangle$. Then $a^2c\in P' \subseteq Q$
and $a^2\in Q$. Thus $c\in Q$. Moreover
$P/\GEN{c}\cong D_8$. By Lemma~\ref{D8}, $a\not\in
Q$. So $(A=\GEN{M,a^2,c},B=\GEN{a^2})$ is a strong
Shoda pair of $H$ with $[H:A]=4$. This  again yields
a contradiction with statement~\ref{Indices}(b) of
Lemma~\ref{SSP-index}.

Now we prove the second statement. Assume that $P$
is  non-abelian and of order $16$. By
statement~\ref{Dihedral} of Lemma~\ref{SSP-index},
$P$ is not isomorphic to $D_{16}$. By the first part
of the lemma, $P$ is not isomorphic to any of the
groups: $Q_{16}$, $D_{16}^-$, $D_{16}^+,D_{16}^-$,
$\mathcal{D}$. The well known description of the
non-abelian groups of order $16$ yields that $P$ is
isomorphic to one of the groups: $Q_8\times C_2$,
$D_8\times C_2$, $\W_{21}$ or $\GEN{a,b \mid
a^{4}=b^{4}=(ab)^{2}=(a^{2},b)=1}$. Hence the result
follows.
\end{proof}

\begin{lemma}\label{NoExp8}
Let $G=G_{k,m,N_2,G_2}$ be a finite group of
Kleinian type. If $G_2$ is non-abelian then its
exponent is $4$, $G_2'\subseteq Z(G_2)$ and $Z(G_2)$
has exponent 2. In particular, $Q_{16}$, $D_{16}^+$,
$D_{16}^-$, $\mathcal{D}$ and $\mathcal{D}^+$ are
not quotients of $G_2$.
\end{lemma}

\begin{proof}
Again we may assume that
$G=G_{0,1,N_2,G_2}=\GEN{w}_3 \rtimes G_2$.

{\it Claim 1}. Let $x,y\in G_{2}$ with $t=(y,x)\neq
1$ and $x$ of order $8$. Then $t$ has order $4$.

In order to prove this we may assume that
$G_{2}=\GEN{x,y}$ and argue by contradiction. So,
suppose that $t$ does not have order $4$. By
statement~\ref{Derivado} of Lemma~\ref{SSP-index}
and statement~\ref{ExpDerivado} of
Lemma~\ref{PropEl}, $t\in Z(G)$ and $t$ has order
$2$. Let $\V=(\GEN{s}_2\times \GEN{y_{1}^2}_4 \times
\GEN{y_{2}^2}_4): (\GEN{\overline{y_1}}_2 \times
\GEN{\overline{y_2}}_2)$, with $s=(y_2,y_1)$ and
$Z(\V)=\GEN{s,y_1^2,y_2^2}$ (this is the same group
$\V$ of Theorem~\ref{Main} with generators renamed
to avoid confusions with the elements $t,x$ and $y$
of $G$). Then, there is an epimorphism $\V
\rightarrow G_2$ mapping $y_1$ to $x$ and $y_2$ to
$y$. Since $\V/\GEN{y_{2}^{2},sy_{1}^{4}}$ has order
$16$ and exponent $8$, $G_2/\GEN{y^2,tx^4}$ has
order at most $16$. However, if
$|G_2/\GEN{y^2,tx^4}|=16$ then $G_2/\GEN{y^2,tx^4}
\cong \V/\GEN{y_2^2,sy_1^4}$ and hence
$G_2/\GEN{y^2,tx^4}$ has exponent $8$, contradicting
Lemma 6.3. This implies that $G_2/\GEN{y^2}$ has
order at most $16$ and $G_2/\GEN{tx^4}$ has order at
most $32$. Since the latter is non-abelian of
exponent $8$, it has order $32$, by
Lemma~\ref{Elementales}.  This implies that
$y^2\not\in\GEN{x,t}$. Indeed, for otherwise
$|G_2|\le 32$ and hence $|G_2 | =
|G_2/\GEN{tx^4}|=32$. So  $t=x^4$ and therefore
$y^2\in \GEN{x}$. Thus $|G_2|=16$, a contradiction.
We thus obtain that $G_{2}/\GEN{y^{2},tx^{4}}$ has
order $8$ because we have seen that this group has
order at most 8 and $G_2/\GEN{tx^4}$ has order 32.
Moreover, since $|G_{2}/\GEN{y^{2}}|\le 16$, using
again Lemma~\ref{Elementales}, the group
$G_{2}/\GEN{y^{2}}$  is either abelian or has
exponent $4$ and thus either $t\in \GEN{y^{2}}$ or
$x^{4}\in \GEN{y^{2}}$. Since $y^{2}\not\in
\GEN{t,x}$, either $t=y^{4}$ or $x^{4}=y^{4}$. So in
both cases we get $x^2y^2\not\in\GEN{tx^4}$ and
$x^4y^4\in\GEN{tx^4}$. This implies that
$G_{2}/\GEN{tx^{4},x^{2}y^{2}}$ has order $16$ and
exponent $8$,  because $x^4\not\in\GEN{tx^4}\cup
\GEN{tx^4}x^2y^2 = \GEN{tx^4,x^2y^2}$.
Lemma~\ref{Elementales} therefore yields that
$G_{2}/\GEN{tx^{4},x^{2}y^{2}}$ is abelian, that is,
$t\in \GEN{tx^{4},x^{2}y^{2}}$. Since $t\not\in
\GEN{tx^4}$ we conclude that $y^2\in\GEN{x,t}$ a
contradiction. This proves the claim.

{\it Claim 2}. If $x\in G_{2}$ has order $8$ then
$(x,(x,G_{2}))=1$.

It is sufficient to show that if $y\in G_{2}$ and
$t=(y,x)\neq 1$ then $(x,t)=1$. Assume the contrary,
then by Lemma~\ref{TodosLosCasosRango2},
$(x,t)=t^2\ne 1$. Hence Claim 1 implies that both
$t$ and $t^2$ have order $4$, a contradiction. This
proves Claim 2.

We now first prove by contradiction that $G_{2}$ has
exponent $4$.  So assume $x\in G_{2}$ has order $8$.
Because of statement~\ref{ExpCentro} of
Lemma~\ref{SSP-index}, we know that $x\not\in
Z(G_{2})$. Let $y\in G_{2}$ be so that $t=(y,x)\neq
1$. As before, we may assume that $G_{2} = \langle
x,y \rangle$.  Because of Claim 1, $t$ has order $4$
and by the second claim $(x,t)=1$.  By
statement~\ref{Derivado} of Lemma~\ref{SSP-index},
$\GEN{t^2}$ is a normal subgroup of $G$ contained in
$N_2$. Then, applying Claim 1 to $G/\GEN{t^2} =
G_{0,1,N_2/\GEN{t^2},G_2/\GEN{t^2}}$, we get that
$x^4\in \GEN{t^2}$. This implies that $t^2=x^4$.
Since $t$ is not central in $\GEN{x,y}$ (as $t$ has
order $4$), we get that $(y,t)\ne 1$ and
$(yx,t)=(yx,(yx,x))\neq 1$. Because of Claim 2 we
obtain that $y^4=(yx)^4=1$. Moreover, by
part~\ref{DosCasos}(a) of Lemma~\ref{PropEl},
$(x^2,y)=t^2$ and, by part~\ref{DosCasos}(b) of the
same lemma, $(y^2,x)=1$. This implies that
$y^{2},tx^{2}\in Z(G_2)$.  Since $t\not\in Z(G_2)$,
we thus have that $G_2/\GEN{y^{2},tx^{2}}$ is a
non-abelian quotient of $D_{16}$. Since $D_{16}$ is
not of Kleinian type, $G_2/\GEN{y^{2},tx^{2}}$ has
order $8$ and, from $t^2=x^4$ and $y^4=1$, we have
that $G_2$ has order at most 32. By
Lemma~\ref{Elementales}, it follows that $G_2$ has
order exactly $32$. Therefore $\GEN{y^2,tx^{2}}$ has
order $4$ and thus $\GEN{y^{2}}\cap \GEN{tx^{2}} =
1$. Hence, both $G_2/\GEN{y^{2}}$ and
$G_2/\GEN{tx^{2}}$ are non-abelian groups of order
$16$. Therefore, by Lemma~\ref{Elementales}, both
have exponent $4$. Thus $x^4\in \GEN{y^{2}}\cap
\GEN{tx^{2}} = 1$, a contradiction. This finishes
the proof of the fact that the exponent of $G$ is
$4$.

We now prove that $G_2'\subseteq Z(G_2)$. We argue
by contradiction. So, because of
statement~\ref{Derivado} in Lemma~\ref{SSP-index},
there exist $x,y\in G_2$ such that $t=(y,x)$ has
order $4$. One may assume without loss of generality
that $x$ and $y$ satisfy condition 1 of
Lemma~\ref{TodosLosCasosRango2} (see
Remark~\ref{t4}), that is $(x,t)=(y,t)=t^2$ and
$x^2,y^2\in Z(G)$. Then $1=(xy)^4 = x t x y^2 x t x
y^2 = t^2$, a contradiction.

It remains to show that $Z(G_2)$ has exponent $2$.
By means of contradiction assume that there exists
$z\in Z(G_2)$ of order 4. Since $G_2$ is not
abelian, there exist $x,y\in G_2$ with $(x,y)=t\neq
1$. As before, one may assume that
$G_2=\GEN{x,y,z}$. Since $t$ has order $2$ and $z$
has order $4$, $H=\GEN{t,z^2,x^2,y^2}$ is an
elementary abelian $2$-subgroup of $Z(G)$ and so
there is a subgroup $L$ of index $2$ in $H$ which
contains $tz^2$ but does not contain $t$. We will
use the bar notation for the natural images of the
elements of $G$ in $G/L$. If $x^2\in L$ we set
$x_1=\overline{x}$, and otherwise we put
$x_1=\overline{tx}$. Similarly, define
$y_1=\overline{y}$ if $y^2\in L$, and
$y_1=\overline{ty}$ otherwise. Then
$G_2/L=\GEN{x_1,y_1,\overline{z}}$ is a non-abelian
epimorphic image of ${\cal D}$ with a central
element $\overline{z}$ of order 4. This yields a
contradiction with Lemma~\ref{Elementales}, because
$L\subseteq N_2$. This finishes the proof.
\end{proof}

We are ready to finish the proof of
Theorem~\ref{Main} by proving that if $G$ is a
non-nilpotent group of Kleinian type then $G$
satisfies condition (F.4). Recall that
$G=G_{k,m,N_2,G_2}$ as in (\ref{NNilp}). Of course,
$G_2$ may be abelian or non-abelian.

Assume first that $G_2$ is abelian. Then
$Z(G)=K\times N_2$. Let $u$ be an element of
minimal order in $G_2\setminus N_2$. Then
$G_2=L\times \GEN{u}$ and $N_2=L\times \GEN{u^2}$.
Because of Lemma~\ref{C_8C_8}, the exponent of
$G_2$ divides 8 and, by statement~\ref{ExpCentro}
of Lemma~\ref{SSP-index}, the exponent of $Z(G)$
divides $4$ or $6$. We separately deal with the
cases $K=1$ and $K\ne 1$. Assume that $K=1$. Then
$G=L\times (M\rtimes \GEN{u})$ and the exponent of
$L$ divides $4$. Therefore $G$ is an epimorphic
image of $A\times H$, where $A$ is abelian of
exponent $4$ and $H$ satisfies the first condition
of (F.4). Assume now that $K\ne 1$, then the
exponent of $Z(G)$ divides $6$ and thus the order
$n$ of $u$ divides $4$. Thus $G$ is an epimorphic
image of $A\times H_1$, with $A$ abelian of
exponent $6$ and $H_1 = M\rtimes C_{n} =
G_{0,m,\GEN{u^2},\GEN{u}_n}$. Then $H_1$ is an
epimorphic image of $H=M\rtimes \W_{11} =
G_{0,m,\GEN{y_1,t,u^2},\W_{11}}$. We conclude that
$G$ is an epimorphic image of $A\times H$, where
$A$ and $H$ satisfy the second condition of (F.4).

Now suppose that $G_2$ is not abelian. Notice that
$Z(G_2)\subset N_2$ because $N_2$ is abelian and
$[G_2:N_2]=2$. By Lemma~\ref{NoExp8}, $G_{2}$ has
exponent $4$ and $G_{2}'$ has exponent $2$.
Furthermore, if $T$ is a proper subgroup of $G_2'$
then $G/T \cong G_{k,m,N_2/T,G_2/T}$ and hence,
{again by Lemma~\ref{NoExp8},} the exponent of
$Z(G_2/T)$ is $2$. It thus follows from
Lemma~\ref{PrimerCaso} that $G_2$ is an epimorphic
image of either $C_2^n\times \W$, $\W_{1n}$ or
$\W_{2n}$ for some $n$.

Assume that $G_2$ is an epimorphic image of
$C_2^n\times \W$, but not of $\W_{in}$ for $i=1,2$
and some $n$. This implies that $G_2=C_2^r\times
\W$ for some $r$. Having in mind that
$Z(G_2)\subseteq N_2$ one has that $G=A\times
G_{0,m,Q,\W}$ for $A$ an abelian group of exponent
dividing $6$ and $Q$ and abelian subgroup of index
2 in $\W$. Let $L_1=\GEN{x^2,y^2},
L_2=\GEN{x^2,(xy)^2}$ and $L_3=\GEN{y^2,(xy)^2}$.
Then $L_i\subseteq Q$ and $\W/L_i\cong D_8$.
Further, the image of $xy$ (resp. $y$, $x$) in
$\W/L_1$ (resp. $\W/L_2$, $\W/L_3$) has order 4.
Thus, $xy,x,y\not\in Q$, contradicting the fact the
$[\W:Q]=2$.

In the remainder of the proof we assume that $G_2$
is an epimorphic image of $\W_{1n}$ or $\W_{2n}$.
For simplicity, the symbols used for the generators
in the description of the groups $\W_{1n}$ and
$\W_{2n}$, as given in part (F) of
Theorem~\ref{Main}, also will be  used for their
images in $G_2$. So we write
$G_2=\GEN{x,y_{1},\cdots ,
y_{n},t_1=(y_1,x),\dots,t_n=(y_n,x)}$ with the
respective relations. Then $G_2'$ is an elementary
abelian $2$-group and $G_2'=\GEN{t_1,\dots,t_n}$.
Assume that $|G_2'|=2^r$. Then, reordering the
$y_i$'s, one may assume that
$G_2'=\GEN{t_1,\dots,t_r}$. Let $r<i\leq n$ and let
$t_i=t_1^{\alpha_1}\cdots t_r^{\alpha_r}$, with
$\alpha_i=0$ or $1$. Then $y'_i=y_i
y_1^{\alpha_1}\cdots y_r^{\alpha_r} \in
Z(G_2)\subseteq N_2$, for $i>r$. Thus, replacing
$y_i$ by $y'_i$ for $i>r$ one has that $G_2 =
B\times P$, where $B=\GEN{y_{r+1},\ldots,y_n}$, an
elementary abelian $2$-group, and $P$ an epimorphic
image of $W_{1r}$ or $W_{2r}$ such that
$P'=\GEN{t_1}_2\times \dots \times \GEN{t_r}_2$.
Then the map $f:P'\rightarrow \GEN{y_1,\ldots,y_r}$
given by $f(t_1^{\alpha_1}\cdots t_r^{\alpha_r}) =
y_1^{\alpha_1}\cdots y_r^{\alpha_r}$ ($\alpha_i=0$
or $1$) is well defined. Moreover $(x,f(s))=s$ and
therefore $(xf(s))^2=sx^2f(s)^2$, for every $s\in
P'$.

Assume that $P$ is a quotient of $\W_{1r}$. Let
$A_1=K\times B$, an abelian group of exponent
dividing $6$. Then $G=A_1\times H_1$, where
$H_1=G_{0,m,Q,P}$  and $Q$ is an abelian subgroup
of index $2$ in $P$. We will show that $G$ is an
epimorphic image of $A\times H$, with $A$ and $H$
satisfying the second condition of $(F.4)$. For
this it is enough to show that one may assume that
$y_1,\dots,y_r,t_1,\dots,t_r,x^2\in Q$. Obviously
$t_1,\dots,t_r,x^2\in Q$. Assume that
$y_{i_0}\not\in Q$. We separately deal with the
cases $x^2\in \GEN{t_{i_0}}$ and $x^2\not\in
\GEN{t_{i_0}}$. If $x^2\not\in \GEN{t_{i_0}}$ then
$K_1=\GEN{x,y_{i_0}}/\GEN{x^2}$ and
$K_2=\GEN{x,y_{i_0}}/\GEN{x^2t_{i_0}}$ are
isomorphic to $D_8$. Moreover $\overline{x}$ has
order $4$ in $K_2$ and $\overline{xy_{i_0}}$ has
order $4$ in $K_1$. Hence $x,xy_{i_0}\not\in Q$, by
Lemma~\ref{D8}, and $y_{i_0}\not\in Q$, yielding a
contradiction. Suppose now that $x^2\in
\GEN{t_{i_0}}$. Then, replacing $x$ by $xt_{i_0}$,
if needed, one may assume that $x^2=1$. Then, for
every $i=1,\dots,r$, we have $\GEN{x,y_i}\cong D_8$
and $xy_i$ has order $4$. Therefore $xy_i\not\in
Q$, by Lemma~\ref{D8}. Since $y_{i_0}\not\in Q$,
one gets that $x\in Q$. For every $1\neq y\in
\GEN{y_1,\dots,y_r}$, the group $\GEN{x,y}$ is not
abelian. Then $y\not \in Q$. That is
$\GEN{y_1,\dots,y_r}\cap Q=1$ and hence $r=1$. Thus
$P=\GEN{x,y_1}\cong D_8$, with $x^2=1=y_i^2$ and
$Q=\GEN{x,t}$. Interchanging the roles of $x$ and
$y_1$  one may assume that $y_1\in Q$ as desired.

Finally, assume that $P$ is a quotient of
$\W_{2r}$, with $|P'|=2^r$. Hence, $f(s)^2 = s$ for
every $s\in P'$. We also assume that $P$ is not an
epimorphic image of $\W_{1h}$ for any $h\ge 1$. We
claim that if $r=1$ then the exponent of $A_1$
divides $2$. Notice that $P$ is non-abelian,
$\W_{21}$ has order $16$ and $D_8$ is an epimorphic
image of $\W_{11}$. Then $P$ is isomorphic to
either $\W_{21}$ or $Q_8$. This implies that $K
\times H_2$ is an epimorphic image of $G$, where
$H_2=G_{0,1,\GEN{a},Q_8}$ and $a$ is an element of
order $4$ in $Q_8$. Then $H_2$ has a cyclic
subgroup $K_2$ of index $2$ and so $(K_2,1)$ is a
strong Shoda pair of $H_2=G_{0,1,\GEN{a},Q_8}$.
Thus $e=e(H_2,K_2,1)$ is a primitive central
idempotent of $\Q H_2$ and, applying
Proposition~\ref{e-pci-meta}, one has $\Q H_2 e
\cong \HQ(\Q(\sqrt{3}))$. Therefore, if $K\neq 1$
then $\Q G$ has a quotient isomorphic to
$\Q(\xi_3)\otimes_{\Q} \HQ(\Q(\sqrt{3})) \simeq
M_2(\xi_3,\sqrt{3})$, contradicting
statement~\ref{Indices} of Lemma~\ref{SSP-index}.
This proves the claim.

Now we separately deal with the cases $x^2\not\in
P'\setminus \{1\}$ and $x^2\in P'\setminus \{1\}$.
In both cases we will show that $r=1$ and hence, by
the above, $K=1$ and $G=B\times H_1$, where $B$ is
an elementary abelian $2$-group and
$H_1=G_{0,m,Q,P}$ with $Q$ is an abelian subgroup
of index $2$ in $P$. Then, in order to show that
$G$ is an epimorphic image of $A\times H$ with $A$
and $H$ satisfying the third condition of (F.4), it
is enough to prove that one may assume that $x\in
Q$. Suppose $x^2\not\in P'\setminus \{1\}$. Then
$\GEN{x,f(s)}/\GEN{x^2}$ is isomorphic to $D_8$,
for every $s\in P'\setminus\{1\}$. By
Lemma~\ref{D8}, $f(s)\not\in Q$ because the natural
image of $f(s)$ in $\GEN{x,f(s)}/\GEN{x^2}$ has
order 4. This implies that $r=1$, because otherwise
$y_1,y_2,y_1y_2\not\in Q$, contradicting the fact
that $[P:Q]=2$. Then, replacing $x$ by $xy_1$ if
needed, one may assume that $x\in Q$. Assume now
that $x^2\in P'\setminus\{1\}$. We claim that one
may assume that $x^2=t_1$. If $x^2 \not \in
\GEN{t_2,\ldots,t_r}$ this is obtained by replacing
$y_1$ by $f(x^2)$. Otherwise $x^2=t_2^{\alpha_2}
\cdots t_r^{\alpha_r}$, for some
$\alpha_1,\ldots,\alpha_r\in \{0,1\}$ with
$\alpha_i=1$ for some $i$. Then, replacing
$\{y_1,\ldots,y_r\}$ by
$\{f(x^2),y_1,\ldots,y_{i-1},y_{i+1},\ldots,y_r\}$,
one obtains the desired conclusion. So we assume
that $x^2=t_1$. Then $\overline{f(s)}$ has order
$4$ in $\GEN{x,f(s)}/\GEN{x^2} \simeq D_8$ for
every $s\in P'\setminus \GEN{t_1}$ and therefore
$(P'\setminus \GEN{t_1})\cap Q = 1$. This implies
that $r\le 2$. If $r=2$ then $y_1y_2,y_2\not\in Q$
and therefore $y_1\in Q$. Replacing $x$ by $xy_2$
if needed, one may assume that $x\in Q$. So
$Q=\GEN{x,y_1,y_2^2}$. Let $1\ne m\in M$. Then
$(U=\GEN{m,x,y_2^2},\GEN{y_2^2})$ is a strong Shoda
pair of $H=\GEN{m,P}$ and $[H:U]=4$, contradicting
statement~\ref{Indices} of Lemma~\ref{SSP-index}.
Thus $r=1$ and $x^2=t_1=y_1^2$. Therefore $P\cong
Q_8$ and either $x$ or $y_1$  does not belong to
$Q$. By symmetry, one may assume that $y_1\not\in
Q$ and, replacing $x$ by $xy_1$ if needed, one may
assume that $x\in Q$. This finishes the proof of
Theorem~\ref{Main}.

\begin{tabular}{ll}
E. Jespers  & A. Pita and A. del R\'{\i}o \\ Dep.
Mathematics &  Dep. Matem\'{a}ticas
\\ Vrije Universiteit Brussel &
Universidad de Murcia \\ Pleinlaan 2 &Campus de
Espinardo\\ 1050 Brussel, Belgium & 30100 Murcia,
Spain\\ efjesper@vub.ac.be &  antopita@um.es,
adelrio@fcu.um.es\\ &\\

 Manuel Ruiz & P. Zalesski\\  Dep. M\'{e}todos
Cuantitativos e Inform\'{a}ticos & Departamento de
Matem\'{a}tica
\\
Universidad Polit\'{e}cnica de Cartagena &
Universidade de Bras\'{\i}lia \\ Paseo Alfonso XIII,
50  &  \\ 30.203 Cartagena, Spain & 70.910-900
Brasilia-DF, Brasil \\ Manuel.Ruiz@upct.es &
pz@mat.unb.br
\end{tabular}

\end{document}